\documentclass[12pt,a4paper]{article}

\usepackage{amsmath}
\usepackage{amssymb}
\usepackage{amsfonts}
\usepackage{graphicx}

\newtheorem{theorem}{Theorem}
\newtheorem{lemma}{Lemma}[section]

\newtheorem{cor}{Corollary}[section]
\newcommand{\Pf}{\noindent\textbf{Proof. }}

\newcommand{\NN}{\mathbb{N}}
\newcommand{\ZZ}{\mathbb{Z}}

\newcommand{\CC}{\mathbb{C}}

\newcommand{\ul}{\underline}

\newcommand\non{\nonumber}
\newcommand\half{\frac{1}{2}}

\def\qed{{\hspace*{\fill}{\vrule height 1ex width 1ex }\quad}
    \vskip 0pt plus20pt}

\def\sgn{{\rm \sgn}}

\begin{document}

\title{A determinantal identity}
\author{Tony C. Dorlas\footnote{Dublin Institute for Advanced Studies, School of Theoretical Physics, 10 Burlington road, Dublin 04, Ireland;} }

\maketitle

\begin{abstract}

We prove an interesting identity for the sum of determinants, which is a generalization of the sum of a geometric progression. 
The proof is quite long and a number of other identities are proved along the way. 
Some of the more elementary ones are deferred to another section at the end.

\end{abstract}



\section{The identity}

We prove an identity for determinants which generalizes the sum of a geometric progression. It was inspired by the problem of calculating scalar products of Bethe Ansatz eigenfunctions, see \cite{Gaudin1981} and \cite{Korepin1982}. 
The proof is quite long and a sequence of lemmas used in the proof is deferred to a second section.
For general properties of determinants, see \cite{Serre2010} and \cite{Greub1967}.

\begin{theorem} For any $n, N \in \NN$ and $a_1,\dots,a_n \in \CC$, the following identity holds.
\begin{eqnarray} && \sum_{1 \leq x_1 < \dots < x_n \leq N} \left| \begin{array}{cccc} a_1^{x_1} & a_1^{x_2} & \cdots &
a_1^{x_n} \\ a_2^{x_1} & a_2^{x_2} & \cdots & a_2^{x_n} \\ \vdots & \vdots & \cdots & \vdots \\ a_n^{x_1} & a_n^{x_2} &
\cdots & a_n^{x_n} \end{array} \right| \non \\ && \qquad = \prod_{k=1}^n \frac{a_k}{a_k-1} \sum_{J \subset
\{1,\dots,n\}}  (-1)^{\nu(J^c)} \gamma(J) \gamma(J^c)  \prod_{j\in J} a_j^N, \label{detid}
\end{eqnarray}
where
\begin{equation} \gamma(J) = \frac{1}{\prod\limits_{\{i,j\} \subset J} (a_i a_j - 1)} \Delta(J), \end{equation}
where \begin{equation} \Delta(J) = \left| \begin{array}{cccc} 1 & a_{j_1} & \cdots & a_{j_1}^{k-1} \\ 1 & a_{j_2}  &
\cdots & a_{j_2}^{k-1} \\ \vdots & \vdots & \cdots & \vdots \\ 1 & a_{j_k} & \cdots & a_{j_k}^{k-1} \end{array} \right|
\end{equation} is a Van der Monde determinant if $J = \{j_1,\dots,j_k\}$ with $j_1 < \dots < j_k$
(and $\gamma(\emptyset) = \gamma(\{k\}) = 1$), and where
\begin{equation} \nu(I) = \sum_{k \in I} k. \end{equation}
\end{theorem}

The proof uses a number of elementary results about determinants of this type, which are stated in
Section~\ref{Lemmasec}.

\Pf We first sum over $x_1$ to write
\begin{eqnarray} && \prod_{k=1}^n \frac{a_k-1}{a_k} \sum_{1 \leq x_1 < \dots < x_n \leq N}
\left| \begin{array}{cccc} a_1^{x_1} & a_1^{x_2} & \cdots & a_1^{x_n} \\
a_2^{x_1} & a_2^{x_2} & \cdots & a_2^{x_n} \\ \vdots & \vdots & \cdots & \vdots \\ a_n^{x_1} & a_n^{x_2} & \cdots &
a_n^{x_n} \end{array} \right| \non \\ &=& \sum_{1 \leq x_2 < \dots < x_n \leq N-1}
\left| \begin{array}{ccccc} a_1^{x_2} - 1 & a_1^{x_2+1} - 1 & a_1^{x_3}(a_1-1) & \cdots & a_1^{x_n}(a_1-1) \\
a_2^{x_2} - 1 & a_2^{x_2+1} - 1 & a_2^{x_3}(a_2-1) & \cdots & a_2^{x_n}(a_2-1) \\
\vdots & \vdots & \vdots & \cdots & \vdots \\ a_n^{x_2} - 1 & a_n^{x_2+1} - 1 & a_n^{x_3}(a_n-1) & \cdots &
a_n^{x_n}(a_n-1) \end{array} \right|. \non \\ &&  \end{eqnarray}

For $n=2$ this becomes
\begin{eqnarray*} && \sum_{x_2=1}^{N-1} \left\{ (a_1 a_2)^{x_2} \left| \begin{array}{cc} 1 & a_1 \\ 1 & a_2 \end{array} \right| -
\left| \begin{array}{cc} 1 & a_1^{x_2+1} - a_1^{x_2} \\ 1 & a_2^{x_2+1} - a_2^{x_2} \end{array} \right| \right\} \\
&=&  \frac{(a_1 a_2)^N - a_1 a_2}{a_1 a_2 - 1} \left| \begin{array}{cc} 1 & a_1 \\ 1 & a_2 \end{array}
\right| - \left| \begin{array}{cc} 1 & a_1^N - a_1 \\ 1 & a_2^N - a_2 \end{array} \right|  \\
&=& \frac{(a_1 a_2)^N - 1}{a_1 a_2 - 1} \left| \begin{array}{cc} 1 & a_1 \\ 1 & a_2 \end{array} \right| + a_1^N - a_2^N
\\ &=& \sum_{J \subset \{1,2\}} (-1)^{\nu(J^c)} \gamma(J) \gamma(J^c) \prod_{j \in J} a_j^N.
\end{eqnarray*}
In general, we want to prove that
\begin{align} &\sum_{1 \leq x_2 < \dots < x_n \leq N-1}
\left| \begin{array}{ccccc} a_1^{x_2} - 1 & a_1^{x_2+1} - 1 & a_1^{x_3}(a_1-1) & \cdots & a_1^{x_n}(a_1-1) \\
a_2^{x_2} - 1 & a_2^{x_2+1} - 1 & a_2^{x_3}(a_2-1) & \cdots & a_2^{x_n}(a_1-1) \\
\vdots & \vdots & \vdots & \cdots & \vdots \\ a_n^{x_2} - 1 & a_n^{x_2+1} - 1 & a_n^{x_3}(a_n-1) & \cdots &
a_n^{x_n}(a_n-1) \end{array} \right| \non \\
&\qquad = \sum_{J \subset \{1,\dots,n\}}  (-1)^{\nu(J^c)} \gamma(J) \gamma(J^c)  \prod_{j\in J} a_j^N. \label{firstsum}
\end{align} We proceed by induction on $n$.
First note that if $x_i=x_{i-1}$, the $i$-th column and the $i-1$th columns are equal (except for $i=3$, in which case
the third column equals the difference of the second and first columns). We can therefore extend the sums to those cases and write the left-hand side of (\ref{firstsum}) as
\begin{eqnarray} &&\sum_{0 \leq x_2 \leq \dots \leq x_n \leq N-1}
\left| \begin{array}{ccccc} a_1^{x_2} - 1 & a_1^{x_2+1} - 1 & a_1^{x_3}(a_1-1) & \cdots & a_1^{x_n}(a_1-1) \\
a_2^{x_2} - 1 & a_2^{x_2+1} - 1 & a_2^{x_3}(a_2-1) & \cdots & a_2^{x_n}(a_2-1) \\
\vdots & \vdots & \vdots & \cdots & \vdots \\ a_n^{x_2} - 1 & a_n^{x_2+1} - 1 & a_n^{x_3}(a_n-1) & \cdots &
a_n^{x_n}(a_n-1) \end{array} \right|\!. \non \\ && \label{gensum} \end{eqnarray}

Expanding the left-hand side of (\ref{gensum}) according to the last column it becomes
\begin{eqnarray} &&\sum_{k=1}^n (-1)^{n-k} (a_k-1) \sum_{x_n=0}^{N-1} a_k^{x_n} \sum_{0 \leq x_2 \leq \dots \leq x_{n-1} \leq x_n} \non \\ &&\quad  \times
\left| \begin{array}{ccccc} a_1^{x_2} - 1 & a_1^{x_2+1} - 1 & a_1^{x_3}(a_1-1) & \cdots & a_1^{x_{n-1}}(a_1-1) \\
a_2^{x_2} - 1 & a_2^{x_2+1} - 1 & a_2^{x_3}(a_2-1) & \cdots & a_2^{x_{n-1}}(a_2-1) \\
\vdots & \vdots & \vdots & \cdots & \vdots \\ \left[ a_k^{x_2} - 1 \right. & a_k^{x_2+1} - 1 & a_k^{x_3}(a_k-1) & \cdots & \left. a_k^{x_{n-1}}(a_k-1) \right] \\
\vdots & \vdots & \vdots & \cdots & \vdots \\ a_n^{x_2} - 1 & a_n^{x_2+1} - 1 & a_n^{x_3}(a_n-1) & \cdots & a_n^{x_{n-1}}(a_n-1) \end{array} \right|, \end{eqnarray} 
where the square brackets around the $k$-th row indicate that
this row is omitted. By the induction hypothesis, this equals
\begin{align} &\sum_{k=1}^n (-1)^{n-k} (a_k-1) \sum_{x_n=0}^{N-1} a_k^{x_n} \non \\ 
&\qquad \times \sum_{J \subset \{1,\dots,n\}\setminus \{k\}} (-1)^{\nu_k(J^c)}\, \gamma(J)\, \gamma(J^c\setminus\{k\})
\prod_{j \in J} a_j^{x_n}, \label{indhyp} \end{align} where $\nu_k(J^c)$ is given by \begin{equation}  \nu_k(I) =
\sum_{i \in I:\,i<k} i + \sum_{i \in I:\, i > k} (i-1). \label{nukdef} \end{equation}

Now let us first consider the case that $|J| = n-1$, that is $J = \{1,\dots,n\} \setminus \{k\}$. In that case
$\nu_k(J^c) = 0$ for all $k$. The corresponding term is
$$ \sum_{k=1}^n (-1)^{n-k} (a_k-1) \sum_{x_n=0}^{N-1} a_k^{x_n}
\gamma(\{1,\dots,n\}\setminus \{k\}) \prod_{j \neq k} a_j^{x_n}. $$ 
Multiplying by $\prod_{1 \leq i < j \leq n} (a_i a_j - 1)$ this becomes
\begin{eqnarray} && \sum_{k=1}^n (-1)^{n-k} (a_k-1) \sum_{x_n=0}^{N-1} a_k^{x_n} \non \\ && \qquad \qquad \times
\prod_{j \neq k} a_j^{x_n} \prod_{j \neq k} (a_j a_k-1)   \, \left| \begin{array}{cccc} 1 & a_1 &
\cdots & a_1^{n-2} \\ \vdots & \vdots & \cdots & \vdots \\ \left[1 \right. & a_k  & \cdots & \left. a_k^{n-1} \right] \\
\vdots & \vdots & \cdots & \vdots \\ 1 & a_n & \cdots & a_n^{n-2} \end{array} \right| \non \\
&=& \sum_{x_n=0}^{N-1} (a_1 \dots a_n)^{x_n} \, \left|
\begin{array}{ccccc} 1 & a_1 & \cdots & a_1^{n-2} & (a_1-1) \prod_{j=2}^n (a_1 a_j-1) \\
\vdots & \vdots & \cdots & \vdots & \vdots \\
1 & a_k  & \cdots & a_k^{n-2} & (a_k-1) \prod_{\substack{j=1 \\ j \neq k}}^n (a_j a_k - 1) \\
\vdots & \vdots & \cdots & \vdots & \vdots \\
1 & a_n & \cdots & a_n^{n-2} & (a_n-1) \prod_{j=1}^{n-1} (a_j a_n-1) \end{array} \right|. \non \\ && \end{eqnarray}
Using Lemma~\ref{L4} and summing over $x_n$ this equals
$$ \left(\prod_{j=1}^n a_j^{N} - 1 \right)
\left| \begin{array}{ccccc} 1 & a_1 & \cdots & a_1^{n-2} & a_1^{n-1} \\
1 & a_2 & \cdots & a_2^{n-2} & a_2^{n-1} \\ \vdots & \vdots & \cdots & \vdots & \vdots \\
1 & a_n & \cdots & a_n^{n-2} & a_n^{n-1} \end{array} \right|. $$ 
The term with $\prod_{j=1}^n a_j^{N}$ is just the
term $|J|=n$ of the right-hand side of (\ref{firstsum}) when divided again by 
$\prod_{1 \leq i < j \leq n}(a_i a_j -1)$. The second term contributes to $J = \emptyset$.

Next consider the case $|J| = n-2$ in the expression (\ref{indhyp}). This equals
\begin{eqnarray*}  &&\sum_{k=1}^n (-1)^{n-k} (a_k-1) \sum_{x_n=0}^{N-1} a_k^{x_n}  
\sum_{l \neq k} (-1)^{\nu_k(\{k,l\})} \gamma(\{k,l\}^c) \prod_{j \neq k,l} a_j^{x_n} \\
&=& \sum_{l=1}^n \sum_{k \neq l} (-1)^{n-k} (a_k-1) \sum_{x_n=0}^{N-1} \prod_{j\neq l} a_j^{x_n} (-1)^{\nu_k(\{k,l\})} \gamma(\{k,l\}^c). \end{eqnarray*} 
Here $\nu_k(J)$ is given by (\ref{nukdef}). Multiplying the $l$-th term by
$\prod\limits_{\substack{1 \leq i < j \leq n \\ i,j \neq l}} (a_i a_j - 1)$ the resulting 
expression is
\begin{align*} &\sum_{k \neq l} (-1)^{n-k} (a_k-1) \sum_{x_n=0}^{N-1} \prod_{j\neq l} a_j^{x_n}
(-1)^{\nu_k(\{k,l\})} \\ &\qquad \times \prod_{\substack{j=1 \\ j\neq k,l}}^n (a_j a_k - 1)\,
\left| \begin{array}{cccc} 1 & a_1 & \cdots & a_1^{n-3} \\ \vdots & \vdots & \cdots & \vdots \\
\left[1 \right. & a_k  & \cdots & \left. a_k^{n-3} \right] \\
\vdots & \vdots & \cdots & \vdots \\
\left[1 \right. & a_l  & \cdots & \left. a_l^{n-3} \right] \\
\vdots & \vdots & \cdots & \vdots \\ 1 & a_n & \cdots & a_n^{n-3} \end{array} \right|. \end{align*} Now,
$$ \nu_k(\{k,l\}) = \begin{cases} l &\text{if $l < k$;} \\ l-1 &\text{if $l > k$.} \end{cases} $$
But, in the case $l > k$, the $l$-th row is below the $k$-th row so the number of rows below 
the $k$-th is only $n-k-1$. Performing the sums over $k$ and $x_n$, we therefore get
\begin{equation*} (-1)^l \frac{\prod_{j \neq l} a_j^N - 1}{\prod_{j \neq l} a_j - 1} \left|
\begin{array}{ccccc} 1 & a_1 & \cdots & a_1^{n-3} & (a_1-1) \prod_{j \neq 1,l} (a_1 a_j -1) \\
\vdots & \vdots & \cdots & \vdots & \vdots \\ \left[1 \right. & a_l  & \cdots & a_l^{n-3} & \left. (a_l-1) \prod_{j \neq l} (a_j a_l-1) \right] \\
\vdots & \vdots & \cdots & \vdots & \vdots \\ 1 & a_n & \cdots & a_n^{n-3} & (a_n-1) \prod_{j \neq l,n} (a_j a_n -1)
\end{array} \right|.
\end{equation*}
We set $J' = J \cup \{k\} = \{l\}^c$ and note that $\nu(J^{\prime c}) = \nu(\{l\}) = l$. Using Lemma~\ref{L4} again, we obtain
$$ \left( \prod_{j \neq l} a_j^N - 1 \right) (-1)^{\nu(\{l\})} \left| \begin{array}{ccccc} 1 & a_1 & \cdots & a_1^{n-3} & a_1^{n-2} \\
\vdots & \vdots & \cdots & \vdots & \vdots \\ \left[1 \right. & a_l  & \cdots & a_l^{n-3} & \left. a_l^{n-2} \right] \\
\vdots & \vdots & \cdots & \vdots & \vdots \\ 1 & a_n & \cdots & a_n^{n-3} & a_n^{n-2} \end{array} \right|.
$$
Dividing again by $\prod\limits_{\substack{1 \leq i < j \leq n \\ i,j \neq l}} (a_i a_j - 1)$ and summing over $l$ this yields
$$ \sum_{l=1}^n (-1)^{\nu(\{l\})} \gamma(\{l\}^c) \left( \prod_{j \neq l} a_j^N - 1 \right). $$
The term with $\prod_{j \neq l} a_j^N$ is just the term $J'=\{l\}^c$ in the right-hand side of (\ref{firstsum}). The other term contributes to the case $J' = \emptyset$.

We now consider the general case  in (\ref{indhyp}). Again, we want to put $J' = J \cup \{k\}$. 
Then $(J')^c = J^c \setminus \{k\}$. Summing over $x_n$ in (\ref{indhyp}) we get
\begin{align*} &\sum_{\substack{J \subset \{1,\dots,n\} \\ J^c \neq \emptyset}} \sum_{k \in J^c}
\left( \frac{a_k^N \prod_{j \in J} a_j^N - 1}{a_k \prod_{j \in J} a_j - 1} \right) \\
&\qquad\qquad \times (-1)^{n-k} (a_k-1) (-1)^{\nu_k(J^c)} \gamma(J) \gamma(J^c \setminus\{k\}). \end{align*} With $J' =
J \cup \{k\}$ this is
\begin{align*} &\sum_{\substack{J' \subset \{1,\dots,n\} \\ J' \neq \emptyset}}
\left(\frac{\prod_{j \in J'} a_j^N - 1}{\prod_{j \in J'} a_j - 1} \right) \\ &\qquad \times \sum_{k \in J'} (-1)^{n-k}
(a_k-1) (-1)^{\nu_k(J^{\prime c} \cup \{k\})} \gamma(J' \setminus \{k\}) \gamma((J')^c).
\end{align*} As in the case $|J|=n-2$, $\nu_k(J^c) = \nu(J^{\prime c}) - p$, where $p$ is the number of $i \in J^c$ with $i > k$,
which compensates for the number of rows below the $k$-th row omitted in the determinant for $J$. Applying
Lemma~\ref{L4}, we therefore obtain
\begin{equation} \sum_{\substack{J' \subset \{1,\dots,n\} \\ J' \neq \emptyset}}
\left(\prod_{j \in J'} a_j^N - 1 \right)  (-1)^{\nu(J^{\prime c})} \gamma(J') \gamma(J^{\prime c}).
\end{equation}
The terms corresponding to $\prod_{j \in J'} a_j^N$ agree with those in the right-hand side of equation (\ref{firstsum}) with $J \neq \emptyset$, so it remains to show that
\begin{equation} - \sum_{\substack{J' \subset \{1,\dots,n\} \\ J' \neq \emptyset}}
(-1)^{\nu(J^{\prime c})} \gamma(J') \gamma(J^{\prime c}) = (-1)^{\nu(\{1,\dots,n\})} \gamma(\{1,\dots, n\}). \end{equation} Equivalently, with $I = J^{\prime c}$,
\begin{equation} \label{zerosum} \sum_{I \subset \{1,\dots,n\}} (-1)^{\nu(I)} \prod_{i \in I} \prod_{j \in I^c} (a_i a_j -1) \Delta(I) \Delta(I^c)  = 0, \end{equation} where
\begin{equation} \Delta(I) = \left| \begin{array}{cccc} 1 & a_{i_1} & \dots & a_{i_1}^{p-1} \\ 1 & a_{i_2} & \dots &
a_{i_2}^{p-1} \\ \vdots & \vdots & \cdots & \vdots \\ 1 & a_{i_p} & \dots & a_{i_p}^{p-1} \end{array} \right| \quad
\mbox{ for } I= \{i_1,\dots, i_p\}. \end{equation}

If $\nu(\{1,\dots,n\})$ is odd, the identity (\ref{zerosum}) is obvious by interchanging $I$ and $I^c$. If $\nu(\{1,\dots,n\})$ is even, the terms $I$ and $I^c$ are equal to each other, so, by symmetry, we can assume that $|I| \leq |I^c|$. 
Then we can expand $\prod_{i \in I} \prod_{j \in I^c} (a_i a_j -1)$ as follows.
\begin{equation} \label{productexpansion}
\prod_{i \in I} \prod_{j \in I^c} (a_i a_j -1) = \sum_{p=0}^{|I|\,|I^c|} \sum_{K \subset I \times I^c: |K| = p}
(-1)^{|I|\,|I^c|-|K|} \prod_{(i,j) \in K} a_i a_j. \end{equation} Set $k = |I|$ so that $|I^c| = n-k$. We can reorder
the points $i \in I$ such that the number $n_i$ of points $(i,j) \in K$ for given $i \in I$, is non-decreasing. Given a
non-decreasing sequence $(n_r)_{r=1}^k$, put $k_m = \#\{r:\,n_r = m\}$. Clearly, $n_r \leq n-k$ so $m \leq n-k$.
Moreover, $\sum_{m=0}^{n-k} k_m = k$. We can then write
\begin{eqnarray} && \sum_{K \subset I \times I^c: |K| = p} \prod_{(i,j) \in K} a_i a_j \non \\ && \qquad = \sum_{\substack{0 \leq n_1 \leq n_2 \leq \dots \leq n_k \leq n-k \\ \sum_{r=1}^k n_r = p}}
\sum_{\substack{(I_m)_{m=0}^k \in \Pi(I):\\ |I_m|=k_m}} \prod_{m=1}^k \prod_{i \in I_m} 
a_i^m \prod_{r=1}^k \sum_{\substack{J_r \subset I^c:\\ |J_r| = n_r}} \prod_{j \in J_r} a_j, \non \\
\end{eqnarray} where $\Pi(I)$ is the set of partitions of $I$. We define
\begin{equation} S_l = \sum_{J \subset I^c:\,|J|=l} \prod_{j \in J} a_j \mbox{ and } S_0=1,\end{equation}
and
\begin{equation} A_{n_1,\dots,n_k} = \sum_{(I_m)_{m=0}^k \in \Pi(I):\,|I_m|=k_m}\, \prod_{m=1}^k \prod_{i \in I_m} a_i^m, \end{equation}
so that
\begin{equation} \label{Kprod} \sum_{K \subset I \times I^c: |K| = p} \prod_{(i,j) \in K} a_i a_j =
\sum_{\substack{0 \leq n_1 \leq n_2 \leq \dots \leq n_k \leq n-k \\ \sum_{r=1}^k n_r = p}} A_{n_1,\dots,n_k}
\prod_{r=1}^k S_{n_r}. \end{equation} Denote $|\ul{n}| = \sum_{r=1}^k n_r$ and
\begin{equation} \NN^k_\uparrow = \{(n_1,\dots,n_k) \in \ZZ^k:\, 0 \leq n_1 \leq \dots \leq n_k \leq n-k \}
\end{equation} and introduce a lexicographic ordering according to
\begin{equation} \ul{n} < \ul{m} \mbox{ if } n_r = m_r \mbox{ for } r > r_0,\ n_{r_0} < m_{r_0}. \end{equation}
We define a map $\phi: \NN^k_{\uparrow} \to \NN^{n-k}_\uparrow$ by
\begin{equation} \phi(\ul{n})_m = \begin{cases} 0 &\text{ if  $1 \leq m \leq n-k-n_k$,} \\
1 &\text{ if $n-k-n_k + 1 \leq m \leq n-k-n_{k-1}$,} \\ \vdots  \\
k &\text{ if $m \geq n-k-n_1 + 1$.} \end{cases} \end{equation} That is, \begin{equation} \phi(\ul{n})_m = \min \{r \geq
0:\, n_{k-r} \leq n-k-m\}. \end{equation} For example, if $n = 10$ and $k=4$, then 
$\phi(0,3,3,5) = (0,1,1,3,3,3)$. A pictorial representation of this map is obtained by filling 
squares of a $k \times (n-k)$ grid with beads; $n_i$ on
column $i$ and $\phi(\ul{n})_j$ on row $j$. We order the rows from bottom to top: 

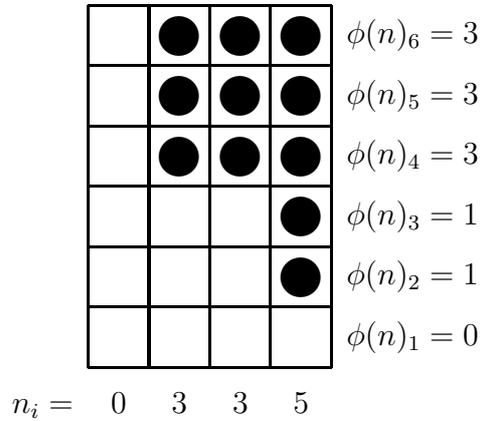
\begin{center}
\setlength{\unitlength}{1mm}
\begin{picture}(60,60) \thicklines
\multiput(10,10)(8,0){5}{\line(0,1){48}} \multiput(10,10)(0,8){7}{\line(1,0){32}} \multiput(38,22)(0,8){5}{\circle*{7}}
\multiput(30,38)(0,8){3}{\circle*{7}} \multiput(22,38)(0,8){3}{\circle*{7}} \put(0,4){$n_i =$} \put(13,4){0}
\put(21,4){3} \put(29,4){3} \put(37,4){5} \put(44,13){$\phi(n)_1 = 0$} \put(44,21){$\phi(n)_2 = 1$}
\put(44,29){$\phi(n)_3 = 1$} \put(44,37){$\phi(n)_4 = 3$} \put(44,45){$\phi(n)_5 = 3$} \put(44,53){$\phi(n)_6 = 3$}
\end{picture}
\end{center}

\textbf{Example 1.} Consider the case where $k=3$, $n=7$ and $I=\{1,2,3\}$. Then (\ref{productexpansion}) reads
\begin{eqnarray*} \lefteqn{\prod_{i \in I} \prod_{j \in I^c} (a_i a_j - 1)} \\ &=& A_{4,4,4} S_4^3  - A_{3,4,4} S_3 S_4^2 \\
&& \quad + A_{3,3,4} S_3^2 S_4 + A_{2,4,4} S_2 S_4^2 \\ && \quad - A_{3,3,3} S_3^3 - A_{2,3,4} S_2 S_3 S_4 - A_{1,4,4}
S_1 S_4^2 \\ && \quad + A_{2,3,3} S_2 S_3^2 + A_{2,2,4} S_2^2 S_4 + A_{1,3,4} S_1 S_3 S_4 + A_{0,4,4} S_4^2 \\
&& \quad -A_{2,2,3} S_2^2 S_3 - A_{1,3,3} S_1 S_3^2 - A_{1,2,4} S_1 S_2 S_4 - A_{0,3,4} S_3 S_4 \\ && \quad + A_{2,2,2}
S_2^3 + A_{1,2,3} S_1 S_2 S_3 + A_{0,3,3} S_3^2 + A_{1,1,4} S_1^2 S_4 + A_{0,2,4} S_2 S_4 \\ && \quad - A_{1,2,2} S_1
S_2^2 - A_{1,1,3} S_1^2 S_3 - A_{0,2,3} S_2 S_3 - A_{0,1,4} S_1 S_4 \\ && \quad + A_{1,1,2} S_1^2 S_2 + A_{0,2,2} S_2^2
+ A_{0,1,3} S_1 S_3 + A_{0,0,4} S_4 \\ && \quad - A_{1,1,1} S_1^3 - A_{0,1,2} S_1 S_2 - A_{0,0,3} S_3 \\ && \quad +
A_{0,1,1} S_1^2 + A_{0,0,2} S_2 - A_{0,0,1} S_1 + 1. \end{eqnarray*} Here we have ordered the terms first according to
$p$ from largest ($p_{\rm max} = k(n-k) = 12$) to smallest ($p=0$) and then according to the above lexicographic
ordering. Here, for example, $A_{0,2,4} = a_1^4(a_2^2 + a_3^2) + a_2^4 (a_1^2 + a_3^2) + a_3^4 (a_1^2 + a_2^2)$, and
$S_2 = \sum_{4 \leq j_1 < j_2 \leq 7} a_{j_1} a_{j_2}$.
\medskip

Let us also define, for $n \in \NN$ and $\ul{m} \in \NN^{n}_0$,
\begin{equation} \label{Deltam} \Delta_{\ul{m}}(J) = \left| \begin{array}{ccc} a_{j_1}^{m_1} & \cdots & a_{j_1}^{m_{n}} \\ \vdots &
\cdots & \vdots \\ a_{j_{n}}^{m_1} & \cdots & a_{j_{n}}^{m_{n}}  \end{array} \right| \mbox{ if } J = \{j_1,\dots,j_n\}.
\end{equation} 
We now claim that the following identities hold.

\begin{lemma} Define $\tilde{m}_r = m_r + r - 1$ and  $\psi(\ul{m}) = \widetilde{\phi(\ul{m})}$,
 i.e. $\psi(\ul{m})_r = \phi(\ul{m})_r + r-1$.
 There is an upper-triangular matrix $R_{\ul{n},\ul{m}}$ such that
\begin{equation} \label{Sproduct} \prod_{r=1}^k S_{n_r} \Delta(I^c) = \sum_{\substack{\ul{m} \in \NN^k_{\uparrow}: \\ |\ul{m}| =
|\ul{n}|}} R_{\ul{n},\ul{m}} \Delta_{\psi(\ul{m})}(I^c), \end{equation}
Moreover, 
\begin{equation} \label{ADelta} A_{\ul{n}} \Delta(I) =
\sum_{\substack{\ul{m} \in \NN^k_\uparrow: \\ \ul{m} \leq \ul{n};\, |\ul{m}| = |\ul{n}|}} (R^{-1})_{\ul{m},\ul{n}} \Delta_{\tilde{\ul{m}}}(I), \end{equation} 
 Equivalently,
\begin{equation} \label{RADel} \sum_{\substack{\ul{n} \in \NN^k_\uparrow: \\ \ul{n} \leq \ul{m};\, |\ul{n}| = |\ul{m}|}}
R_{\ul{n},\ul{m}} A_{\ul{n}} \Delta(I) =  \Delta_{\tilde{\ul{m}}}(I). \end{equation}
\end{lemma}

\Pf
The identity (\ref{Sproduct}) follows by induction from Corollary~\ref{L6cor} of Lemma~\ref{L6}, according to which
\begin{equation} \label{Slformula}
S_l \Delta_{\ul{m}}(I^c) = \sum_{\substack{\ul{m}': \, m'_i-m_i = 0,1 \\ 
\sum_{i=1}^{n-k} (m'_i-m_i) = l}} \Delta_{\ul{m}'}(I^c) \quad (\ul{m} \in \NN_\uparrow^{n-k}).\end{equation}
(Note that we have to replace $k$ by $n-k$ in Corollary~\ref{L6cor}.)
In terms of the pictorial representation, this means that multiplication by $S_l$ corresponds to the addition of $l$ additional beads on the right-most empty sites of $l$ different rows such that the resulting sequence is still non-decreasing.  For example, in the tableau for $\phi(0,3,3,5)$ with $n=10$ above, in case $l=3$, we could add the three new beads shown as unfilled circles:

\begin{center}
\setlength{\unitlength}{1mm}
\begin{picture}(60,60) \thicklines
\multiput(10,10)(8,0){5}{\line(0,1){48}} \multiput(10,10)(0,8){7}{\line(1,0){32}} 
\multiput(38,22)(0,8){5}{\circle*{7}}
\multiput(30,38)(0,8){3}{\circle*{7}} \multiput(22,38)(0,8){3}{\circle*{7}} 
\multiput(30,22)(0,8){2}{\circle{5}} \put(14,54){\circle{5}}
\put(0,4){$n'_i =$} \put(13,4){1} \put(21,4){3} \put(29,4){5} \put(37,4){5} 
\put(44,13){$\phi(n')_1 = 0$} \put(44,21){$\phi(n')_2 = 2$}
\put(44,29){$\phi(n')_3 = 2$} \put(44,37){$\phi(n')_4 = 3$} \put(44,45){$\phi(n')_5 = 3$} \put(44,53){$\phi(n')_6 = 4$}
\end{picture}
\end{center}

This is equivalent to adding a total of $l$ beads on the upper most empty sites of a number of columns  such that there are no new beads horizontally next to each other. Thus we can also write
\begin{equation} \label{Sop} S_l \Delta_{\psi(\ul{m})}(I^c) =
\sum_{\substack{\ul{m}': \, m_i \leq m'_i \leq m_{i+1} \\ \sum_{i=1}^k (m'_i-m_i) = l}} \Delta_{\psi(\ul{m}')}(I^c);\quad (\ul{m} \in \NN_\uparrow^k).
\end{equation}
In particular, note that the minimal $\ul{m}'$ (w.r.t. the above ordering) is obtained by adding beads to the upper-most incomplete rows. Note also that $|\ul{m}'| = |\ul{m}| + l$. 

Iterating this formula, it follows that for $\ul{n} \in \NN^k_{\uparrow}$,
\begin{equation} \label{Sprod} \prod_{r=1}^k S_{n_r} \Delta(I^c) = \sum_{\substack{\ul{m} \in \NN^k_\uparrow: \\ \ul{m} \geq \ul{n};\,|\ul{m}| = |\ul{n}|}}
R_{\ul{n},\ul{m}} \Delta_{\psi(\ul{m})} (I^c), \end{equation} where the matrix $R$ is upper-triangular and has integer
matrix elements given by the number of times a given configuration $\ul{m}$ is obtained by iterating the above
procedure. (Note that the number of non-zero $n_r$ is the maximal length of a row of beads, i.e. $\phi(\ul{n})_{n-k}$.
Also, $\ul{m} = \ul{n}$ only if the beads are placed in order from right to left starting with $n_k$, so
$R_{\ul{n},\ul{n}} = 1$.)
\medskip

\textbf{Example 2.} In the case of Example 1, with $p=6$, the matrix $R$ is given by
$$ R = \left( \begin{array}{ccccc} 1 & 2 & 1 & 1 & 3 \\ 0 & 1 & 1 & 1 & 2 \\ 0 & 0 & 1 & 0 & 1 \\ 0 & 0 & 0 & 1 & 1 \\ 0 & 0 & 0 & 0 & 1 \end{array} \right), $$ 
where the rows are numbered from top to bottom and the columns from left to
right in increasing lexicographic order, i.e. $(222), (123), (033), (114), (024)$.
For example, consider the matrix element $R_{(123),(024)}$ (right-most element of the second row). 
Applying $S_3$ we obtain three beads in the upper squares of the last column of a $3 \times 4$ tableau.
Then applying $S_2$ we obtain two different tableaus: one with two beads on the upper squares 
of the second column in addition to the original 3 beads, and one with one bead on the upper square
of the second column and 4 beads on the third column. Both these tableaus are less in lexicographic order than $(024)$. Applying next $S_1$ there is only one way to obtain $(024)$ in each case. 
Therefore $R_{(123),(024)} = 2$.
\medskip

Next we consider the expressions $A_{\ul{n}} \Delta(I)$. We prove (\ref{ADelta}) by induction on
$p=|\ul{m}|$ and $k$. 
(Note that for $p=1$, we have $\ul{n} = (0,\dots,0,1) = \ul{m}$ and $A_{0,\dots,0,1} =
\sum_{i=1}^k a_i$ so that $A_{0,\dots,0,1} \Delta(I) = \Delta_{0,1,\dots,k-2, k}$ follows from
Corollary~\ref{L5cor} of Lemma~\ref{L5}. If $k=1$ then $I = \{i\}$ and $A_{n_1} = a_i^{n_1}$ so $A_{n_1} \Delta(\{i\}) = A_{n_1} = \Delta_{n_1}(\{i\})$.)

If $m_1 \geq 1$ then we define $\ul{m}''$ by $m''_i = m_i - m_1$. Now it is easy to see that \begin{equation}
R_{\ul{n},\ul{m}} = 0 \mbox{ if } \sum_{r \leq r_0} n_r < \sum_{r \leq r_0} m_r \mbox{ for some } r_0\geq 1.
\end{equation} In particular, if $m_1 \geq 1$ then $n_1 \geq m_1$, and in that case
\begin{equation} R_{\ul{n},\ul{m}} = R_{\ul{n}'',\ul{m}''}, \end{equation} where $n''_i = n_i - m_1$.
Since $|\ul{n}''| = |\ul{n}|-k m_1$, it follows from the induction hypothesis (w.r.t. $p$) that
\begin{eqnarray} \sum_{\substack{\ul{n} \in \NN^k_\uparrow: \\ \ul{n} \leq \ul{m};\, |\ul{n}| = \ul{m}|}}
R_{\ul{n},\ul{m}} A_{\ul{n}} \Delta(I) &=& A_{m_1,\dots,m_1} \sum_{\substack{\ul{n}'' \in \NN^k_\uparrow: \\ \ul{n}''
\leq \ul{m}'';\, |\ul{n''}| = |\ul{m}|-k m_1}} R_{\ul{n}'',\ul{m}''} A_{\ul{n}''} \Delta(I) \non \\ &=&
A_{m_1,\dots,m_1} \Delta_{\ul{m}''} = \Delta_{\ul{m}}. \end{eqnarray}

It therefore remains to consider the case that $m_1=0$.  If $m_1 = 0$ and also $n_1=0$, then we can define $\ul{m}' = (m_2,\dots,m_k)$ and $\ul{n}' = (n_2,\dots,n_k)$ so that 
\begin{equation} R^{(k)}_{\ul{n},\ul{m}} = R^{(k-1)}_{\ul{n}',\ul{m}'}. \end{equation}
We can therefore also assume that $n_1 \geq 1$.

For ease of notation, we can assume that $I = \{1,\dots,k\}$. By induction with respect to $k$ and expanding with respect to the first column, it follows that
\begin{eqnarray} && \sum_{\substack{\ul{n} \in \NN^k_\uparrow: \\ \ul{n} \leq \ul{m};\, \sum_{r=1}^k n_r = p}}
R_{\ul{n},\ul{m}} A_{\ul{n}} \Delta(I) = \Delta_{\tilde{\ul{m}}}(I) \non \\ && \qquad +
\sum_{i=1}^k (-1)^{i-1} \sum_{s=1}^{m_k} a_i^s \sum_{\substack{\ul{n} \in \NN^k_\uparrow: \\
\ul{n} \leq \ul{m};\, |\ul{n}| = |\ul{m}|,\,s \in \{n_i\}}}  R_{\ul{n},\ul{m}} A_{\ul{n}^{(s)}}
\Delta'(I\setminus\{i\}),
\end{eqnarray}
where $\Delta'(I\setminus\{i\})$ denotes $\Delta_{1,2,\dots,k-1}(I\setminus\{i\})$, and $\ul{n}^{(s)}$ is obtained from
$\ul{n}$ by omitting $s$, i.e. if $n_r = s$ then $n^{(s)}_i = n_i$ for $i<r$ and $n^{(s)}_i = n_{i+1}$ if $i \geq r$.
(There may be more than one such $r$, namely, if $k_s > 1$, in which case we can simply choose one.) By the definition
of $R$, and the formula (\ref{Sop}) for $S_l$, we have that if $s \in \{n_i\}$ then
\begin{equation} \label{Riteration} R_{\ul{n},\ul{m}} = \sum_{\substack{\ul{q} \in \NN^{k-1}_{\uparrow}: \\ |\ul{q}| =
s, \, 0 \leq q_i \leq m_{i+1}-m_i}} R_{\ul{n}^{(s)},\ul{m}'-\ul{q}}. \end{equation} Inserting this, the remainder term
becomes
\begin{eqnarray} && \sum_{i=1}^k (-1)^{i-1} \sum_{s=1}^{m_k} a_i^s 
\sum_{\substack{\ul{n} \in \NN^k_\uparrow: \\
\ul{n} \leq \ul{m};\, |\ul{n}| = |\ul{m}|,\,s \in \{n_i\}}}  R_{\ul{n},\ul{m}} A_{\ul{n}^{(s)}} \Delta'(I\setminus\{i\}) \non \\
&& = \sum_{s=1}^{m_k} \sum_{\substack{\ul{q} \in \NN^{k-1}_{\uparrow}: \\ |\ul{q}| = s, \, 0 \leq q_i \leq
m_{i+1}-m_i}} \sum_{i=1}^k (-1)^{i-1}  a_i^s \non \\ && \qquad\qquad \times \sum_{\substack{\ul{n}' \in \NN^{k-1}_\uparrow: \\
\ul{n}' \leq \ul{m}'-\ul{q};\, |\ul{n}'| = |\ul{m}|-s}}  R_{\ul{n}',\ul{m}'-\ul{q}} A_{\ul{n}'}
\Delta'(I\setminus\{i\}) \non \\ && = \sum_{s=1}^{m_k} \sum_{\substack{\ul{q} \in \NN^{k-1}_{\uparrow}: \\ |\ul{q}| =
s, \, 0 \leq q_i \leq m_{i+1}-m_i}} \sum_{i=1}^k (-1)^{i-1}  a_i^s
\Delta'_{\widetilde{\ul{m}'-\ul{q}}}(I\setminus\{i\}) \non
\\ && = \sum_{s=1}^{m_k} \sum_{\substack{\ul{q} \in
\NN^{k-1}_{\uparrow}: \\ |\ul{q}| = s, \, 0 \leq q_i \leq m_{i+1}-m_i}} 
\Delta_{\widetilde{(s,\ul{m}'-\ul{q})}}(I).
\label{remDeltas} \end{eqnarray}
(The second equality follows from the induction hypothesis.)
\medskip

\textbf{Example 3.} To clarify this, consider Example 1 again and let $\ul{m} = (024)$. With the $R$-matrix of Example
2, we then have
\begin{equation*} \sum_{\substack{\ul{n} \in \NN^3_\uparrow: \\ \ul{n} \leq \ul{m};\, |\ul{n}| = |\ul{m}|}}
R_{\ul{n},\ul{m}} A_{\ul{n}} \Delta(\{1,2,3\}) = \left| \begin{array}{ccc} \sum_{\ul{n} \leq (024);\,|\ul{n}| = 6}
R_{\ul{n},024} A_{\ul{n}} & a_1 & a_1^2 \\ \sum_{\ul{n} \leq (024);\,|\ul{n}| = 6} R_{\ul{n},024} A_{\ul{n}} & a_2 &
a_2^2 \\ \sum_{\ul{n} \leq (024);\,|\ul{n}| = 6} R_{\ul{n},024} A_{\ul{n}} & a_3 & a_3^2 \end{array} \right|.
\end{equation*} In the $i$-th row, we separate out the terms  where $a_i$ has the power 0 in $A_{\ul{n}}$, in
particular $n_1=0$. These are given by $A_{\ul{n}'}$ (as a function of $a_j\,(j \neq i)$). This yields
\begin{equation*} \left|
\begin{array}{cccc} \sum_{\ul{n}' \leq (24);\,|\ul{n}'| = 6} R_{\ul{n}',24} A_{\ul{n}'} & a_1 & a_1^2 \\ \sum_{\ul{n}'
\leq (24);\,|\ul{n}| = 6} R_{\ul{n}',24} A_{\ul{n}'} & a_2 & a_2^2 \\ \sum_{\ul{n}' \leq (24);\,|\ul{n}| = 6}
R_{\ul{n}',24} A_{\ul{n}'} & a_3 & a_3^2  \end{array} \right| = \left|
\begin{array}{cccc} 1 & a_1^3 & a_1^6 \\ 1 & a_2^3 & a_2^6 \\ 1 & a_3^3 & a_3^6  \end{array} \right| =
\Delta_{\tilde{\ul{m}}}(\{1,2,3\}). \end{equation*} 
(For example, $$ \sum_{\ul{n}' \leq 24;\, |\ul{n}| = 6} R_{\ul{n}',24} A_{\ul{n}'}(a_2,a_3) 
\left| \begin{array}{cc} a_2 & a_2^2 \\ a_3 & a_3^2 \end{array} \right| =
 \left| \begin{array}{cc} a_2^3 & a_2^6 \\ a_3^3 & a_3^6 \end{array} \right|. $$
 by the induction hypothesis.)

The remaining terms have $a_i^s$ for some $s\geq 1$ in the first column of the $i$-th row. They are
$$ \sum_{s=1}^4 a_i^s \sum_{\substack{\ul{n} \in \NN^3_\uparrow:\\ \ul{n} \leq (024);\,|\ul{n}| = 6,\,s \in \{n_i\}}}
R_{\ul{n},024} A_{\ul{n}}. $$ Now, by equation (\ref{Riteration}),
\begin{eqnarray*} R_{222,024} &=& R_{22,04} + R_{22,13} + R_{22,22} = 3; \\
R_{123,024} &=& R_{23,14} + R_{23,23} = R_{13,04} + R_{13,13} = R_{12,03} + R_{12,12} = 2; \\
R_{033,024} &=& R_{33,24} = R_{03,03} = 1; \\
R_{114,024} &=& R_{14,14} = R_{11,02} =1 \mbox{ and } R_{024,024} = R_{24,24} = R_{04,04} = R_{02,02} = 1.
\end{eqnarray*}
(Note that $R_{03,12} = 0$ for example, and in the case of $R_{114,024}$, the term $R_{11,11}$ is not allowed because
in that case $q_2=3$ whereas $m_3-m_2=2$.)

For $s=4$ we thus obtain $a_i^4 (A_{02} + R_{114,024} A_{11}) = a_i^4 (A_{02} + R_{11,02} A_{11})$, which yields the
determinant
$$ \left| \begin{array}{ccc} a_1^4 (A_{02} + R_{11,02} A_{11}) & a_1 & a_1^2 \\ a_2^4 (A_{02} + R_{11,02} A_{11}) &
a_2 & a_2^2 \\ a_3^4 (A_{02} + R_{11,02} A_{11}) & a_3 & a_3^2 \end{array} \right| = \Delta_{4,1,4} = 0. $$ (Here we
use induction w.r.t. $k$.) For $s=3$ we obtain $a_i^3 (A_{03} + R_{123,024} A_{12}) = a_i^3 (A_{03} + R_{12,03} A_{12}
+ R_{12,12} A_{12})$. This yields the determinant
$$ \left| \begin{array}{ccc} a_1^3 (A_{03} + R_{12,03} A_{12} + R_{12,12} A_{12}) & a_1 & a_1^2 \\
a_2^3 (A_{03} + R_{12,03} A_{12} + R_{12,12} A_{12}) & a_2 & a_2^2 \\ a_3^3 (A_{03} + R_{12,03} A_{12} + R_{12,12}
A_{12}) & a_3 & a_3^2 \end{array} \right| $$ which equals $\Delta_{3,1,5} + \Delta_{3,2,4} = -\Delta_{1,3,5}
-\Delta_{2,3,4}$. For $s=2$ we get $a_i^2 (A_{04} + R_{123,024} A_{13} + R_{222,024} A_{22}) = a_i^2 (A_{04} +
(R_{13,04} + R_{13,13}) A_{13} + (R_{22,04} + R_{22,13} + R_{22,22}) A_{22})$. This yields $\Delta_{2,1,6} +
\Delta_{2,2,5} + \Delta_{2,3,4} = -\Delta_{1,2,6} + \Delta_{2,3,4}$. Finally, for $s=1$ we have $a_i (R_{114,024}
A_{14} + R_{123,024} A_{23}) = a_i (A_{14} + (R_{23,14} + R_{23,23}) A_{23})$ and we obtain the determinants
$\Delta_{1,2,6} + \Delta_{1,3,5}$. In total, we get $-\Delta_{1,3,5} -\Delta_{2,3,4} -\Delta_{1,2,6} + \Delta_{2,3,4} + \Delta_{1,2,6} + \Delta_{1,3,5} = 0$.
This completes the analysis of this example.
\bigskip

In general, we shall prove that the resulting determinants in (\ref{remDeltas}) cancel in pairs. 
Consider a term $\Delta_{\widetilde{(s,\ul{m}'-\ul{q})}}$. 
It equals $\pm \Delta_{\widetilde{\ul{n}}}$ for some $\ul{n} \in \NN^k_\uparrow$. 
Conversely, now first suppose that $\ul{n}' = \ul{m}'-\ul{q}$ for some  $\ul{q}$ satisfying $0 \leq
q_i \leq m_{i+1} - m_i$ and $|\ul{q}| = n_1 \geq 1$, i.e. there is a term with $s=n_1$. 
Then consider the case $s=n_2+1$, where we need  
$\ul{m}' -\tilde{\ul{q}} = (n_1-1,n_3,\dots,n_k)$. Set $\tilde{q}_1 = q_1 + n_2 - n_1 +1$ and
$\tilde{q}_i = q_i$ for $i \geq 2$. Then $(n_1-1,n_3,\dots,n_k) = \ul{m}'-\tilde{\ul{q}}$ and $|\tilde{\ul{q}}| = n_2+1$. Moreover, since $1 \leq n_1 \leq n_2$ and $n_2 = m_2-q_1 \leq m_2$, 
we have $\tilde{q}_1 = m_2-n_1+1 \geq 0$ and $\tilde{q}_1 \leq m_2 = m_2-m_1$. 
It follows that if $\Delta_{\tilde{\ul{n}}}$ occurs in the sum (\ref{remDeltas})
(with $s = n_1$) then $\Delta_{\widetilde{(n_2+1,n_1-1,n_3,\dots,n_k)}}$ also occurs. But they cancel one another.
Conversely, suppose that $\Delta_{\widetilde{(n_2+1,n_1-1,n_3,\dots,n_k)}}$ occurs, so that 
$(n_1-1,n_3,\dots,n_k) = \ul{m}'-\tilde{\ul{q}}$ for some $\tilde{\ul{q}}$ such that 
$0 \leq \tilde{q}_i \leq m_{i+1}-m_i$ and $|\tilde{\ul{q}}|
= n_2+1$. Define $q_1=\tilde{q}_1-n_2+n_1-1$. Then we need that $0 \leq q_1 \leq m_2$, i.e. $0 \leq m_2-n_2\leq m_2$. 
Therefore, if $\Delta_{\widetilde{(n_2+1,n_1-1,n_3,\dots,n_k)}}$ occurs then 
$\Delta_{\tilde{\ul{n}}}$ also occurs provided $n_2 \leq m_2$. 
If this is not the case then $\Delta_{\tilde{\ul{n}}}$ does not occur and we must start
with $\Delta_{\widetilde{(n_2+1,n_1-1,n_3,\dots,n_k)}}$. 
If this term does occur in (\ref{remDeltas}) then there is $\ul{q}$ such that 
$0 \leq q_i \leq m_{i+1}-m_i$,  $|\ul{q}| = n_2+1$ and $(n_1-1,n_3,\dots,n_k) = \ul{m}'-\ul{q}$.
Defining $\tilde{\ul{q}}$ by $\tilde{q}_2 = q_2 + n_3 - n_2 + 1$, $\tilde{q}_i = q_i$ for $i \neq 2$, we have $|\ul{q}| = n_3+2$ and $(n_1-1,n_2-1,n_4,\dots,n_k) = \ul{m}'-\tilde{\ul{q}}$. 
We need $0 \leq \tilde{q}_2 \leq m_3-m_2$, i.e.
$m_2+1 \leq n_2 \leq m_3+1$. But $ n_2 \geq m_2+1$ because we assumed that 
$\Delta_{\tilde{\ul{n}}}$ does not occur. On the other hand $n_2 \leq n_3 = m_3-q_2 \leq m_3$. Therefore, $\Delta_{\widetilde{(n_2+1,n_1-1,n_3,\dots,n_k)}}$ also
occurs, and the two terms cancel each other.

More generally, suppose that $r \geq 2$ is an integer such that $\pm \Delta_{\tilde{\ul{n}}}$ with $s=n_r+r-1$ occurs in the
sum (\ref{remDeltas}). Then there exists $\ul{q} \in \NN_0^{k-1}$ such that $0 \leq q_i \leq m_{i+1}-m_i$ and $|\ul{q}| = n_r + r-1$ and 
$(n_1-1,\dots,n_{r-1}-1,n_{r+1},\dots,n_k) = (m_2-q_1,\dots,m_k-q_{k-1})$. Therefore,
\begin{equation} \label{Delexist} \begin{cases} m_i+1 \leq n_i \leq m_{i+1}+1 &\text{for $i < r$; } \\ m_i \leq n_{i+1} \leq m_{i+1}
&\text{for $i \geq r$.} \end{cases} \end{equation} Define $\ul{q}^{(r-1)}$ by $q^{(r-1)}_{r-1} = q_{r-1} - n_r + n_{r-1} - 1$ and $q^{(r-1)}_i = q_i$ for $i \neq r-1$. 
Then $|\ul{q}^{(r-1)}| = |\ul{q}| - n_r + n_{r-1} + 1 = n_r + r -1 - n_r + n_{r-1} - 1 = 
n_{r-1} + r - 2$ and $m_r-q^{(r-1)}_{r-1} = m_r - q_{r-1} + n_r - n_{r-1} + 1 = n_r$.
Therefore the term with $s=n_{r-1} + r - 2$ also occurs provided $0 \leq q^{(r-1)}_{r-1} \leq m_r - m_{r-1}$. But, $m_r - q_{r-1} = n_r$, so this holds if $0 \leq m_r - n_r \leq m_r - m_{r-1}$. 
By (\ref{Delexist}), $n_{r-1} \geq m_{r-1}+1$
and since $n_r \geq n_{r-1}$ the second inequality holds. 
Thus the term $s = n_{r-1} + r - 2$ also occurs if $n_r \leq m_r$.
\medskip

Suppose now that this term does not occur. Then we conclude that $n_r \geq m_r + 1$.  
Define $\ul{q}^{(r+1)}$ by
$q^{(r+1)}_{r} = q_r + n_{r+1} - n_r + 1$ and $q^{(r+1)}_i = q_i$ for $i \neq r$. 
Then $(n_1-1,\dots,n_r-1, n_{r+2},\dots,n_k) = \ul{m}'-\tilde{\ul{q}}$ since 
$m_{r+1} - q^{(r+1)}_r = m_{r+1} - (q_r + n_{r+1} - n_r + 1) = n_r -1$. 
Also, $|\ul{q}^{(r+1)}| = n_{r+1} + r$. Moreover, $n_r \geq m_r + 1 \implies q^{(r+1)}_r \leq m_{r+1} - m_r$ and
$n_{r+1} \leq m_{r+1} \implies n_r \leq m_{r+1} + 1 \implies q^{(r+1)}_r \geq 0$. Therefore the term with $s = n_{r+1} + r$ also occurs and cancels the term $s = n_r + r -1$.

We conclude that if the term $s = n_r + r - 1$ occurs then either $s = n_{r-1} + r - 2$ exists or 
$s = n_{r+1} + r $ exists, but not both. Note that $s \leq m_k$, so only terms $s = n_r + r - 1$ can exist where $n_r \leq m_k$. That
means that if the term $s = n_k + k - 1$ occurs then the term $s = n_{k-1} + k - 2$ also occurs.
This proves that the sum (\ref{remDeltas}) equals zero, and hence that (\ref{RADel}) holds. \qed
\medskip

Now, inserting
(\ref{Sprod}) and (\ref{RADel}) into (\ref{Kprod}) we have
\begin{eqnarray} \lefteqn{\sum_{K \subset I \times I^c: |K| = p} \prod_{(i,j) \in K} a_i a_j \Delta(I) \Delta(I^c) } \non \\
&=& \sum_{\substack{\ul{n} \in \NN^k_\uparrow: \\ |\ul{n}| = p}} A_{\ul{n}} \Delta(I) \sum_{\substack{\ul{m} \geq
\ul{n}: \\ |\ul{m}| = p}} R_{\ul{n},\ul{m}} \Delta_{\psi(\ul{m})}(I^c) \non \\
&=& \sum_{\substack{\ul{m} \in \NN^k_\uparrow: \\ |\ul{m}| = p}} \Delta_{\psi(\ul{m})}(I^c) \sum_{\substack{\ul{n} \leq
\ul{m}: \\ |\ul{n}| = p}} R_{\ul{n},\ul{m}}A_{\ul{n}} \Delta(I) \non \\
&=& \sum_{\substack{\ul{m} \in \NN^k_\uparrow: \\ |\ul{m}| = p}} \Delta_{\tilde{\ul{m}}}(I) \Delta_{\psi(\ul{m})}(I^c).
\end{eqnarray}
Inserting this into (\ref{productexpansion}) and (\ref{zerosum}) we have
\begin{eqnarray} \lefteqn{ \sum_{I \subset\{1,\dots,n\}} (-1)^{\nu(I)} \prod_{i \in I} \prod_{j \in I^c} (a_i a_j -
1) \Delta(I) \Delta(I^c) } \non \\ &=& \sum_{I \subset\{1,\dots,n\}} (-1)^{\nu(I)} \sum_{p=0}^{|I|\,|I^c|}
(-1)^{|I|\,|I^c| - p} \sum_{\substack{\ul{m} \in \NN^k_\uparrow: \\ |\ul{m}| = p}} \Delta_{\tilde{\ul{m}}}(I)
\Delta_{\psi(\ul{m})}(I^c) \non \\ &=& \sum_{k=0}^n \sum_{p=0}^{k(n-k)} (-1)^{k(n-k)-p} \sum_{\ul{m} \in
\NN^k_{\uparrow}; \,|\ul{m}| = p} \non \\ && \qquad \times \sum_{\substack{I \subset \{1,\dots,n\}: \\ |I| = k}}
(-1)^{\nu(I)} \Delta_{\tilde{\ul{m}}}(I) \Delta_{\psi(\ul{m})}(I^c). \end{eqnarray} 
The last sum is an expansion of
$\Delta_{\tilde{\ul{m}}, \psi(\ul{m})}(\{1,\dots,n\})$ with respect to the first $k$ columns. 
In general,
\begin{eqnarray} \lefteqn{\Delta_{m_1,\dots,m_n}(\{1,\dots,n\}) =} \non \\ && = (-1)^{k(k+1)/2} \sum_{\substack{I \subset \{1,\dots,n\}:\\|I|
= k}} (-1)^{\nu(I)} \Delta_{m_1,\dots,m_k}(I) \Delta_{m_{k+1},\dots,m_n}(I^c). \end{eqnarray} Indeed, for $k=1$ we have
$$ \Delta_{m_1,\dots,m_n}(\{1,\dots,n\}) = \sum_{i \in \{1,\dots,n\}} (-1)^{i-1} a_i^{m_1}
\Delta_{m_{2},\dots,m_n}(\{1,\dots,n\}\setminus \{i\}), $$ where $a_i^{m_1} = \Delta_{m_1}(\{i\})$. By induction this yields
\begin{eqnarray*} \Delta_{m_1,\dots,m_n}(\{1,\dots,n\}) &=& \sum_{i=1}^n (-1)^{i-1} a_i^{m_1} \Delta_{m_2,\dots,m_n} (\{1,\dots,n\} \setminus \{i\})
\\ &=& (-1)^{k(k-1)/2} \sum_{i=1}^n (-1)^{i-1} a_i^{m_1} \\
&& \times \sum_{\substack{I \subset \{1,\dots,n\}\setminus\{i\}:\\|I| = k-1}} (-1)^{\nu'_i(I)}
\Delta_{m_2,\dots,m_k}(I) \Delta_{m_{k+1},\dots,m_n}(I^c), \end{eqnarray*} 
where $\nu'_i(I) = \sum_{j \in I} j - \#\{j \in I:\, j>i\}$. Thus
\begin{eqnarray*} \lefteqn{\Delta_{m_1,\dots,m_n}(\{1,\dots,n\}) =} \\ &=& (-1)^{k(k-1)/2}  \sum_{\substack{I \subset \{1,\dots,n\}:\\|I| = k}}
\sum_{i \in I} (-1)^{\nu(I)-1-\#\{j \in I:\,j>i\}} \\ && \qquad \times a_i^{m_1} \Delta_{m_2,\dots,m_k}(I) \Delta_{m_{k+1},\dots,m_n}(I^c)
\\ &=& (-1)^{k(k-1)/2}  \sum_{\substack{I \subset \{1,\dots,n\}:\\|I| = k}}
\sum_{i \in I} (-1)^{\nu(I)-k+\#\{j \in I:\,j<i\}} \\ && 
\qquad \times a_i^{m_1} \Delta_{m_2,\dots,m_k}(I) \Delta_{m_{k+1},\dots,m_n}(I^c)
\\ &=& (-1)^{k(k+1)/2}  \sum_{\substack{I \subset \{1,\dots,n\}:\\|I| = k}}
(-1)^{\nu(I)} \Delta_{m_1,\dots,m_k}(I) \Delta_{m_{k+1},\dots,m_n}(I^c). \end{eqnarray*} 
Hence, in order to prove (\ref{zerosum}), we want to show that
\begin{equation} \label{reducedexpr} \sum_{k=0}^n \sum_{p=0}^{k(n-k)} (-1)^{k(n-k)-p+k(k+1)/2} \sum_{\substack{\ul{m} \in \NN^k_{\uparrow}:
\\|\ul{m}| = p}} \Delta_{\tilde{\ul{m}},\psi(\ul{m})}(\{1,\dots,n\}) = 0. \end{equation}
First note that $\Delta_{\tilde{\ul{m}},\psi(\ul{m})}(\{1,\dots,n\}) = 0$ unless $\tilde{\ul{m}}$ and $\psi(\ul{m})$
have nothing in common and make up $\{0,1,\dots,n-1\}$. In particular, $|\ul{m}| + |\phi(\ul{m})| + \half k(k-1) +
\half (n-k)(n-k-1) = \half n(n-1)$, i.e. $$ 2p = \half n(n-1) - \half k(k-1) - \half (n-k)(n-k-1) = k(n-k). $$ If $k(n-k)$ is odd, there is no nonzero term, so if $n$ is even then $k$ must also be even. 
We therefore need
\begin{equation} \sum_{k=0}^n (-1)^{k(n-k)/2+k(k+1)/2} \sum_{\substack{\ul{m} \in \NN^k_{\uparrow}:
\\|\ul{m}| = k(n-k)/2}} \Delta_{\tilde{\ul{m}},\psi(\ul{m})}(\{1,\dots,n\}) = 0. \end{equation}

Next we argue that $\Delta_{\tilde{\ul{m}},\psi(\ul{m})}(\{1,\dots,n\}) = 0$ unless $m_i + m_{k-i+1} = n-k$ for
$i=1,\dots,k$. Consider the case $i=1$. In order that all the numbers below $m_1$ are present, we need $\phi(\ul{m})_j
= 0$ for $j=1,\dots,m_1$, while $\phi(\ul{m})_{m_1+1} \geq 1$. This means that the number of zeros in $\phi(\ul{m})$
equals $m_1$, so $m_k = n-k-m_1$. The converse also holds. Similarly, for $i>1$, we must have $\psi(\ul{m})_j = j +
i-1$, i.e. $\phi(\ul{m})_j = i$, for $j = m_{i-1}+1,\dots,m_i$ and $\phi(\ul{m})_{m_i+1} \geq i+1$. This implies that
$m_{k-i+2} - m_{k-i+1} = m_i-m_{i-1}$. By induction, therefore $m_i + m_{k-i+1} = m_{i-1} + m_{k-i+2} = n-k$. In
particular, if $k$ is odd, then $n-k$ is even and $m_{(k+1)/2} = (n-k)/2$.

Consider first the case that $n$ is even, and hence also $k$ is even. 
Then we can count the number of possible solutions as follows. We choose the values of $\tilde{m}_i$ with $i=1,\dots,n/2$ arbitrarily between 1 and $n/2$. These are strictly increasing and determine 
uniquely $m_1,\dots,m_{k/2}$. The remaining $m_i$ ($i=k/2+1,\dots,k$) are then given by the condition $m_i + m_{k-i+1} = n-k$, and the values of $\psi(\ul{m})_j$ are given by the interstices. 
The number of possible solutions is therefore $\displaystyle{{n/2 \choose k/2}}$. 

Note also, that if we move the $k/2$ last elements $m_i$ ($i=k/2+1,\dots,k$) across all $\phi(\ul{m})_j$ ($j=1,\dots,n-k$), then in order to put the
$\tilde{m}_i$ and $\psi(\ul{m})_j$ in increasing order, it remains to move each $m_i$ with 
$i \leq k/2$ across equally many $\phi(\ul{m})_j$ to the right as we need to move $m_{k-i+1}$ 
across $\phi(\ul{m})_j$ to the left. This means that in each case, the determinant 
$\Delta_{\tilde{\ul{m}},\psi(\ul{m})} = (-1)^{k(n-k)/2} \Delta(\{1,\dots,n\})$. 
Inserting this into the left-hand side of (\ref{reducedexpr}) we obtain
\begin{eqnarray} && \sum_{k=0}^n (-1)^{k(n-k)/2+k(k+1)/2} 
\sum_{\substack{\ul{m} \in \NN^k_{\uparrow}: \\|\ul{m}| = k(n-k)/2}} 
\Delta_{\tilde{\ul{m}},\psi(\ul{m})}(\{1,\dots,n\}) = \non \\ && \qquad =
\sum_{\substack{k=0 \\ k\,{\rm even}}}^n (-1)^{k/2} {n/2 \choose k/2} \Delta(\{1,\dots,n\}) = 0. \end{eqnarray} (Note
that if $k$ is even ,then $(-1)^{k(k+1)/2} = (-1)^{k/2}$.)

Analogously, if $n$ is odd, then if $k$ is even, the number of possibilities is $\displaystyle{{(n-1)/2 \choose k/2}}$,
and if $k$ is odd then the number of possibilities is $\displaystyle{{(n-1)/2 \choose (k-1)/2}}$. The sign is again
$(-1)^{k(n-k)/2}$ and we obtain \begin{eqnarray} \lefteqn{\sum_{k=0}^{(n-1)/2} (-1)^{k(k+1)/2} {(n-1)/2 \choose [k/2]}
\Delta} \non \\  &=& \sum_{l=0}^{(n-1)/2} (-1)^l {(n-1)/2 \choose l} \Delta + \sum_{l=0}^{(n-1)/2} (-1)^{l+1} {(n-1)/2
\choose l} \Delta = 0. \end{eqnarray} In both cases therefore (\ref{reducedexpr}) holds. The claim (\ref{zerosum}) is
thus proved. This completes the proof of the theorem. \qed

\section{Lemmas} \label{Lemmasec}

\begin{lemma} \label{L1} Let $\cal R$ be a commutative ring. For $n \geq 3$ and $a_1,\dots,a_n \in {\cal R}$, and for
$0 \leq k+l \leq n-2$,
\begin{equation} \left| \begin{array}{ccccc} 1 & a_1 & \cdots & a_1^{n-2} & a_1^l \sum_{2 \leq j_1 < \dots < j_k \leq
n} a_{j_1} \dots a_{j_k} \\ 1 & a_2 & \cdots & a_2^{n-2} & a_2^l \sum_{\substack{1 \leq j_1 < \dots < j_k \leq n \\ j_r
\neq 2}} a_{j_1} \dots a_{j_k} \\ \vdots & \vdots & \cdots & \vdots & \vdots \\ 1 & a_n & \cdots & a_n^{n-2} & a_n^l
\sum_{1 \leq j_1 < \dots < j_k \leq n-1} a_{j_1} \dots a_{j_k} \end{array} \right| = 0. \end{equation}
\end{lemma}

\Pf For $k=0$ this is obvious.

We now proceed by induction on $k$:
\begin{eqnarray*} && \left| \begin{array}{ccccc} 1 & a_1 & \cdots & a_1^{n-2} & a_1^l 
\sum_{2 \leq j_1 < \dots < j_k \leq n} a_{j_1} \dots a_{j_k} \\ 
1 & a_2 & \cdots & a_2^{n-2} & 
a_2^l \sum_{\substack{1 \leq j_1 < \dots < j_k \leq n \\  j_r \neq 2}} a_{j_1} \dots a_{j_k} \\ \vdots & \vdots & \cdots & \vdots & \vdots \\ 1 & a_n & \cdots & a_n^{n-2} & a_n^l
\sum_{1 \leq j_1 < \dots < j_k \leq n-1} a_{j_1} \dots a_{j_k} \end{array} \right| \\ 
&& \qquad = \left|
\begin{array}{ccccc} 1 & a_1 & \cdots & a_1^{n-2} & a_1^l \sum_{1 \leq j_1 < \dots < j_k \leq n} a_{j_1} \dots a_{j_k} \\ 1 & a_2 & \cdots & a_2^{n-2} & a_2^l 
\sum_{1 \leq j_1 < \dots < j_k \leq n} a_{j_1} \dots a_{j_k} \\ \vdots & \vdots
& \cdots & \vdots & \vdots \\ 1 & a_n & \cdots & a_n^{n-2} & a_n^l 
\sum_{1 \leq j_1 < \dots < j_k \leq n} a_{j_1} \dots a_{j_k}
\end{array} \right| \\ && \quad - \left| \begin{array}{ccccc} 1 & a_1 & \cdots & a_1^{n-2} & a_1^{l+1}
\sum_{2 \leq j_1 < \dots < j_{k-1} \leq n} a_{j_1} \dots a_{j_{k-1}} \\ 1 & a_2 & \cdots & a_2^{n-2} & a_2^{l+1}
\sum_{\substack{1 \leq j_1 < \dots < j_{k-1} \leq n \\ j_r \neq 2}} a_{j_1} \dots a_{j_{k-1}} \\ \vdots & \vdots &
\cdots & \vdots & \vdots \\ 1 & a_n & \cdots & a_n^{n-2} & a_n^{l+1} \sum_{1 \leq j_1 < \dots < j_{k-1} \leq n-1}
a_{j_1} \dots a_{j_{k-1}} \end{array} \right| = 0 \end{eqnarray*} provided $k+l \leq n-2$. Indeed, the first term
equals zero because the last column is a constant multiple of the $l+1$-th column, where $l \leq n-2$. The second term
equals zero by the induction hypothesis. \qed

Similarly, we have also
\begin{lemma} \label{L2}
Let $\cal R$ be a commutative ring. For $n \geq 3$ and $a_1,\dots,a_n \in {\cal R}$, and for $1 \leq k,l \leq n-1$,
such that $k+l \geq n$,
\begin{equation} \left| \begin{array}{ccccc} 1 & a_1 & \cdots & a_1^{n-2} & a_1^l \sum_{2 \leq j_1 < \dots < j_k \leq
n} a_{j_1} \dots a_{j_k} \\ 1 & a_2 & \cdots & a_2^{n-2} & a_2^l \sum_{\substack{1 \leq j_1 < \dots < j_k \leq n \\ j_r
\neq 2}} a_{j_1} \dots a_{j_k} \\ \vdots & \vdots & \cdots & \vdots & \vdots \\ 1 & a_n & \cdots & a_n^{n-2} & a_n^l
\sum_{1 \leq j_1 < \dots < j_k \leq n-1} a_{j_1} \dots a_{j_k} \end{array} \right| = 0. \end{equation}
\end{lemma}

\Pf For $l \geq 1$ and $k=n-1$ the final element in the $i$-th row equals $a_i^{l-1} a_1 \dots a_n$ so the determinant
is zero. For $k < n-1$ we write
\begin{eqnarray*} a_i^l \sum_{\substack{1 \leq j_1 < \dots < j_k \leq n \\ j_r \neq i}} a_{j_1} \dots a_{j_k} &=&
a_i^{l-1} \sum_{\substack{1 \leq j_1 < \dots < j_k \leq n \\ j_r \neq i}} a_i a_{j_1} \dots a_{j_k} \\
&=& a_i^{l-1} \sum_{1 \leq j_1 < \dots < j_{k+1} \leq n} a_{j_1} \dots a_{j_{k+1}} \\ && \qquad - a_i^{l-1} \sum_{\substack{1
\leq j_1 < \dots < j_{k+1} \leq n \\ j_r \neq i}} a_{j_1} \dots a_{j_{k+1}}. \end{eqnarray*}
Then first terms inserted into the determinant yield zero since $l-1 \leq n-2$, and the second terms yield zero by induction
provided $l \geq 1$. \qed

\begin{lemma} \label{L3}
Let $\cal R$ be a commutative ring. For $n \geq 3$ and $a_1,\dots,a_n \in {\cal R}$, and for $ 0 \leq k \leq n-2$,
\begin{eqnarray} \lefteqn{\left| \begin{array}{ccccc} 1 & a_1 & \cdots & a_1^{n-2} & a_1^{n-1-k}
\sum_{2 \leq j_1 < \dots < j_k \leq n} a_{j_1} \dots a_{j_k} \\ 1 & a_2 & \cdots & a_2^{n-2} & a_2^{n-1-k}
\sum_{\substack{1 \leq j_1 < \dots < j_k \leq n \\ j_r \neq 2}} a_{j_1} \dots a_{j_k} \\ \vdots & \vdots & \cdots &
\vdots & \vdots \\ 1 & a_n & \cdots & a_n^{n-2} & a_n^{n-1-k}
\sum_{1 \leq j_1 < \dots < j_k \leq n-1} a_{j_1} \dots a_{j_k} \end{array} \right|} \non \\
&& = (-1)^k \left| \begin{array}{ccccc} 1 & a_1 & \cdots & a_1^{n-2} & a_1^{n-1} \\
1 & a_2 & \cdots & a_2^{n-2} & a_2^{n-1} \\ \vdots & \vdots & \cdots & \vdots & \vdots \\
1 & a_n & \cdots & a_n^{n-2} & a_n^{n-1} \end{array} \right|. \end{eqnarray}
\end{lemma}

\Pf For $k=0$ the identity is tautological. For $k\geq 1$ we write 
\begin{eqnarray*} \lefteqn{a_i^{n-1-k} 
\sum_{\substack{1 \leq j_1 < \dots < j_k \leq n \\ j_r \neq i}} a_{j_1} \dots a_{j_k}} \\ &=&
a_i^{n-k-1} \sum_{1 \leq j_1 < \dots < j_{k} \leq n} a_{j_1} \dots a_{j_{k}} - a_i^{n-k}
\sum_{\substack{1 \leq j_1 < \dots < j_{k-1} \leq n \\ j_r \neq i}} a_{j_1} \dots a_{j_{k-1}}. \end{eqnarray*} 
The first term yields zero and the result follows by induction. \qed

As a corollary we have

\begin{lemma} \label{L4} For $n \geq 3$ and $a_1,\dots,a_n \in {\cal R}$,
\begin{eqnarray} \lefteqn{\left| \begin{array}{ccccc} 1 & a_1 & \cdots & a_1^{n-2} & (a_1-1)\prod_{j=2}^n (a_1 a_j - 1)  \\
1 & a_2 & \cdots & a_2^{n-2} & (a_2-1)\prod_{\substack{j=1 \\ j \neq 2}}^n (a_2 a_j - 1) \\
\vdots & \vdots & \cdots & \vdots & \vdots \\
1 & a_n & \cdots & a_n^{n-2} & (a_n-1) \sum_{j=1}^{n-1} (a_j a_n - 1) \end{array} \right| } \non \\
&& = \left(\prod_{i=1}^n a_i - 1\right) \left| \begin{array}{ccccc} 1 & a_1 & \cdots & a_1^{n-2} & a_1^{n-1} \\
1 & a_2 & \cdots & a_2^{n-2} & a_2^{n-1} \\
\vdots & \vdots & \cdots & \vdots & \vdots \\ 1 & a_n & \cdots & a_n^{n-2} & a_n^{n-1} \end{array} \right|. \end{eqnarray}
\end{lemma}

\Pf We expand
$$ \prod_{\substack{j=1 \\ j \neq i}}^n (a_i a_j - 1) = \sum_{k=0}^{n-1} (-1)^{n-k-1} a_i^k \sum_{\substack{1 \leq j_1 < \dots < j_k \leq n \\ j_r \neq i}} a_{j_1} \dots a_{j_k}. $$ 

First consider the case that  $n$ is even.  Consider the term $a_i$ in the factor $a_i-1$.
Then the total power of $a_i$ is $l=k+1$ so $k+l \geq n$ if $k \geq n/2$. By Lemma~\ref{L1} 
and Lemma~\ref{L2} these terms yield zero unless $2k+1=n-1$ or $l=n$ and $k=n-1$. 
The latter is the highest-order term and yields
$$ \prod_{i=1}^n a_i \left| \begin{array}{ccccc} 1 & a_1 & \cdots & a_1^{n-2} & a_1^{n-1} \\ 1 & a_2 & \cdots & a_2^{n-2} & a_2^{n-1} \\
\vdots & \vdots & \cdots & \vdots & \vdots \\ 1 & a_n & \cdots & a_n^{n-2} & a_n^{n-1} \end{array} \right|. $$ 
If $k=n/2-1$, then by Lemma~\ref{L3} this yields the contribution
$$ (-1)^{n-1-k} (-1)^k \left| \begin{array}{ccccc} 1 & a_1 & \cdots & a_1^{n-2} & a_1^{n-1} \\ 1 & a_2 & \cdots & a_2^{n-2} & a_2^{n-1} \\
\vdots & \vdots & \cdots & \vdots & \vdots \\ 1 & a_n & \cdots & a_n^{n-2} & a_n^{n-1} \end{array} \right| =
-\left| \begin{array}{ccccc} 1 & a_1 & \cdots & a_1^{n-2} & a_1^{n-1} \\ 1 & a_2 & \cdots & a_2^{n-2} & a_2^{n-1} \\
\vdots & \vdots & \cdots & \vdots & \vdots \\ 1 & a_n & \cdots & a_n^{n-2} & a_n^{n-1} \end{array} \right|\!. $$
The term -1 in the factor $a_i-1$ does not contribute because $2k \neq n-1$ and $l \leq n-1$.

Next consider the case that $n$ is odd. Then $2k+1 \neq n-1$ so the $a_i$ term only contributes the highest-order term. Setting $l=k=(n-1)/2$ we obtain by Lemma~\ref{L3}, 
$$ -(-1)^{n-1} \left| \begin{array}{ccccc} 1 & a_1 & \cdots & a_1^{n-2} & a_1^{n-1} \\ 1 & a_2 & \cdots & a_2^{n-2} & a_2^{n-1} \\
\vdots & \vdots & \cdots & \vdots & \vdots \\ 1 & a_n & \cdots & a_n^{n-2} & a_n^{n-1} \end{array} \right| =
-\left| \begin{array}{ccccc} 1 & a_1 & \cdots & a_1^{n-2} & a_1^{n-1} \\ 1 & a_2 & \cdots & a_2^{n-2} & a_2^{n-1} \\
\vdots & \vdots & \cdots & \vdots & \vdots \\ 1 & a_n & \cdots & a_n^{n-2} & a_n^{n-1} \end{array} \right|\!. $$
(The minus sign in $a_i-1$ compensates for the fact that $n-1$ is now even.) \qed

\begin{lemma} \label{L5} For $n \geq 2$, $1 \leq k \leq n-1$, and $a_1,\dots,a_n \in {\cal R}$,
\begin{eqnarray} && \left| \begin{array}{cccc} \sum_{2 \leq j_1 < \dots < j_k \leq
n} a_{j_1} \dots a_{j_k} & a_1 & \cdots & a_1^{n-1} \\ \sum_{\substack{1 \leq j_1 < \dots < j_k \leq n \\ j_r \neq
2}} a_{j_1} \dots a_{j_k} & a_2 & \cdots & a_2^{n-1}  \\ \vdots & \vdots & \cdots & \vdots \\
\sum_{\substack{1 \leq j_1 < \dots < j_k \leq n \\ j_r \neq n}} a_{j_1} \dots a_{j_k} & a_n & \cdots & a_n^{n-1}
\end{array} \right| \non \\
&& = \left| \begin{array}{ccccccc} 1 & a_1 & \cdots & a_1^{n-k-1} & a_1^{n-k+1} & \cdots & a_1^n \\
1 & a_2 & \cdots & a_2^{n-k-1} & a_2^{n-k+1} & \cdots & a_2^n \\
\vdots & \vdots & \cdots & \vdots & \vdots & \cdots & \vdots \\
1 & a_n & \cdots & a_n^{n-k-1} & a_n^{n-k+1} & \cdots & a_n^n \end{array} \right|.
\end{eqnarray}
\end{lemma}

\Pf For $k=n-1$, we have, expanding,
\begin{eqnarray*} \lefteqn{\left| \begin{array}{cccc} a_{2} \dots a_{n} & a_1 & \cdots & a_1^{n-1} \\ a_1 a_3 \dots a_n & a_2 & \cdots & a_2^{n-1}  \\
\vdots & \vdots & \cdots & \vdots \\ a_1 \dots a_{n-1} & a_n & \cdots & a_n^{n-1} \end{array} \right|} \\ &=&
\sum_{j=1}^n (-1)^{j-1} a_1 \dots a_{j-1} a_{j+1} \dots a_n \, \left| \begin{array}{ccc}  a_1 & \cdots & a_1^{n-1} \\
\vdots & \cdots & \vdots \\  \left[ a_j \right. & \cdots & \left. a_j^{n-2} \right]  \\ \vdots & \cdots & \vdots \\ a_n
& \cdots & a_n^{n-1} \end{array} \right| \\ &=&
 \sum_{j=1}^n (-1)^{j-1} \left| \begin{array}{ccc}  a_1^2 & \cdots & a_1^{n} \\
\vdots & \cdots & \vdots \\  \left[a_j^2 \right. & \cdots & \left. a_j^{n}\right]  \\ \vdots & \cdots & \vdots \\
a_n^2 & \cdots & a_n^{n} \end{array} \right| = \left| \begin{array}{cccc}  1& a_1^2 & \cdots & a_1^{n} \\
\vdots & \vdots & \cdots & \vdots  \\ 1 & a_n^2 & \cdots & a_n^{n} \end{array} \right|. \end{eqnarray*} 
We proceed by induction and write similarly,
\begin{eqnarray*} \lefteqn{\left| \begin{array}{cccc} \sum_{\substack{1\leq j_1 < \dots < j_k \leq n \\ j_r \neq 1}}
a_{j_1} \dots a_{j_k} & a_1 & \cdots & a_1^{n-1} \\ \sum_{\substack{1\leq j_1 < \dots < j_k \leq n \\ j_r \neq 2}}
a_{j_1} \dots a_{j_k} & a_2 & \cdots & a_2^{n-1}  \\
\vdots & \vdots & \cdots & \vdots \\ \sum_{\substack{1\leq j_1 < \dots < j_k \leq n \\ j_r \neq n}} a_{j_1} \dots
a_{j_k} & a_n & \cdots & a_n^{n-1} \end{array} \right|} \\ &=& \sum_{j=1}^n (-1)^{j-1} \sum_{\substack{1\leq j_1 <
\dots < j_k \leq n \\ j_r \neq j}} a_{j_1} \dots a_{j_k} \, \left| \begin{array}{ccc}  a_1 & \cdots & a_1^{n-1} \\
\vdots & \cdots & \vdots \\  \left[ a_j \right. & \cdots & \left. a_j^{n-1} \right]  \\ \vdots & \cdots & \vdots \\ a_n
& \cdots & a_n^{n-1} \end{array} \right| \\ &=& \sum_{j=1}^n (-1)^{j-1} 
\prod_{i\neq j} a_i \left| \begin{array}{cccc}
\sum_{\substack{1\leq j_1 < \dots < j_k \leq n \\ j_r \neq j}}
a_{j_1} \dots a_{j_k} & a_1 & \cdots & a_1^{n-2} \\
\vdots & \cdots & \vdots \\  \big[ \qquad"  & a_j & \cdots & a_j^{n-2}\big]  \\ 
\vdots & \vdots & \cdots & \vdots \\
\sum_{\substack{1\leq j_1 < \dots < j_k \leq n \\ j_r \neq j}} a_{j_1} \dots a_{j_k} & a_n & \cdots & a_n^{n-2} \end{array} \right| \non \\ &=& 
\sum_{j=1}^n (-1)^{j-1} \prod_{i\neq j} a_i^2 \left| \begin{array}{ccccc}
\sum_{\substack{1\leq j_1 < \dots < j_{k-1} \leq n \\ j_r \neq 1,j}}
a_{j_1} \dots a_{j_{k-1}} & 1 & a_1 & \cdots & a_1^{n-3} \\ \vdots & \cdots & \vdots \\
\bigg[\, \sum_{\substack{1\leq j_1 < \dots < j_k \leq n \\ j_r \neq j}}
a_{j_1} \dots a_{j_k}  & 1 & a_j & \cdots &  a_j^{n-3}\bigg]  \\ 
\vdots & \vdots & \vdots & \cdots & \vdots \\
\sum_{\substack{1\leq j_1 < \dots < j_{k-1} \leq n \\ j_r \neq j,n}} a_{j_1} \dots a_{j_{k-1}} 
& 1 & a_n & \cdots & a_n^{n-3} \end{array} \right| \non \\ && 
+ \sum_{j=1}^n (-1)^{j-1} \prod_{i\neq j} a_i \left| \begin{array}{cccc}
\sum_{\substack{2\leq j_1 < \dots < j_k \leq n \\ j_r \neq j}} a_{j_1} \dots a_{j_k} & a_1 & \cdots & a_1^{n-2} \\ \vdots & \cdots & \vdots \\  
\big[\qquad "  & a_j & \cdots & a_j^{n-2}\big]  \\ \vdots & \vdots & \cdots & \vdots \\
\sum_{\substack{1\leq j_1 < \dots < j_k \leq n-1 \\ j_r \neq j}} a_{j_1} \dots a_{j_k} 
& a_n & \cdots & a_n^{n-2} \end{array} \right|. \end{eqnarray*}
The first term equals zero by Lemma~\ref{L1} since $k-1 \leq n-3$. By the induction hypothesis, the second term equals
\begin{eqnarray*} && \sum_{j=1}^n (-1)^{j-1} \prod_{i\neq j} a_i \left| \begin{array}{ccccccc}
1 & a_1 & \cdots & a_1^{n-k-2} & a_1^{n-k} & \cdots & a_1^{n-1} \\
1 & a_2 & \cdots & a_2^{n-k-2} & a_2^{n-k} & \cdots & a_2^{n-1} \\
\vdots & \vdots & \cdots & \vdots & \vdots & \cdots & \vdots \\ \left[ 1 \right. & a_j & \cdots & a_j^{n-k-2} &
a_n^{n-k} & \cdots & \left. a_j^{n-1} \right] \\ \vdots & \vdots & \cdots & \vdots & \vdots & \cdots & \vdots \\
1 & a_n & \cdots & a_n^{n-k-2} & a_n^{n-k} & \cdots & a_n^{n-1} \end{array} \right| \\ && = \left|
\begin{array}{ccccccc} 1 & a_1 & \cdots & a_1^{n-k-1} & a_1^{n-k+1} & \cdots & a_1^n \\
1 & a_2 & \cdots & a_2^{n-k-1} & a_2^{n-k+1} & \cdots & a_2^n \\
\vdots & \vdots & \cdots & \vdots & \vdots & \cdots & \vdots \\
1 & a_n & \cdots & a_n^{n-k-1} & a_n^{n-k+1} & \cdots & a_n^n \end{array} \right|.
\end{eqnarray*}
\qed

\begin{cor} \label{L5cor} For $n \geq 2$, $1 \leq k \leq n-1$, and $a_1,\dots,a_n \in {\cal R}$,
\begin{eqnarray} && \sum_{1 \leq j_1 < \dots < j_k \leq n} a_{j_1} \dots a_{j_k} \left| \begin{array}{cccc}
1 & a_1 & \cdots & a_1^{n-1} \\ 1 & a_2 & \cdots & a_2^{n-1}  \\ \vdots & \vdots & \cdots & \vdots \\
1 & a_n & \cdots & a_n^{n-1} \end{array} \right| \non \\
&& = \left| \begin{array}{ccccccc} 1 & a_1 & \cdots & a_1^{n-k-1} & a_1^{n-k+1} & \cdots & a_1^n \\
1 & a_2 & \cdots & a_2^{n-k-1} & a_2^{n-k+1} & \cdots & a_2^n \\
\vdots & \vdots & \cdots & \vdots & \vdots & \cdots & \vdots \\
1 & a_n & \cdots & a_n^{n-k-1} & a_n^{n-k+1} & \cdots & a_n^n \end{array} \right|.
\end{eqnarray}
\end{cor}

\Pf Writing $$ \sum_{1 \leq j_1 < \dots < j_k \leq n} a_{j_1} \dots a_{j_k} = a_j \sum_{\substack{1 \leq j_1 < \dots < j_{k-1} \leq n \\ j_r \neq j}} a_{j_1} \dots a_{j_{k-1}} + 
\sum_{\substack{1 \leq j_1 < \dots < j_{k} \leq n \\ j_r \neq j}} a_{j_1} \dots a_{j_k}, $$ 
we see that the first term vanishes if $k \leq n-1$ by Lemma~\ref{L1}. \qed

We generalise Lemma~\ref{L5} further:
\begin{lemma} \label{L6} Let $\cal R$ be a commutative ring and $a_1,\dots,a_n \in {\cal R}$. Let $n \in \NN$ and $m_1,\dots, m_{n-1} \in \NN_0$ such that $1 \leq m_1 < \dots < m_{n-1}$. Then, for any $k \in \NN$ with $1 \leq k \leq n-1$,
\begin{eqnarray} && \left| \begin{array}{cccc} \sum_{2 \leq j_1 < \dots < j_k \leq n} a_{j_1} \dots a_{j_k} & a_1^{m_1} & \cdots & a_1^{m_{n-1}} \\
\sum_{\substack{1 \leq j_1 < \dots < j_k \leq n \\ j_r \neq 2}} a_{j_1} \dots a_{j_k} & a_2^{m_1} & \cdots & a_2^{m_{n-1}}  \\
\vdots & \vdots & \cdots & \vdots \\ \sum_{\substack{1 \leq j_1 < \dots < j_k \leq n \\ j_r \neq n}} a_{j_1} \dots
a_{j_k} & a_n^{m_1} & \cdots & a_n^{m_{n-1}} \end{array} \right| \non \\
&& = \sum_{\substack{m_1 \leq m'_1 < \dots < m'_{n-1}: \, (\forall i) m'_i - m_i = 0,1 \\  \#\{i:\,m'_i = m_i + 1\} =
k}} \left| \begin{array}{cccc} 1 & a_1^{m'_1}  & \cdots & a_1^{m'_{n-1}} \\
\vdots & \vdots & \cdots & \vdots \\ 1 & a_n^{m'_1} & \cdots & a_n^{m'_{n-1}}  \end{array} \right|.
\end{eqnarray}
\end{lemma}

\Pf We proceed as in the previous lemma and first note that
\begin{eqnarray*} && \left| \begin{array}{cccc} \sum_{2 \leq j_1 < \dots < j_{n-1} \leq n} a_{j_1} \dots a_{j_k} & a_1^{m_1} & \cdots & a_1^{m_{n-1}} \\
\sum_{\substack{1 \leq j_1 < \dots < j_{n-1} \leq n \\ j_r \neq 2}} a_{j_1} \dots a_{j_k} & a_2^{m_1} & \cdots & a_2^{m_{n-1}}  \\
\vdots & \vdots & \cdots & \vdots \\ 
\sum_{\substack{1 \leq j_1 < \dots < j_{n-1} \leq n \\ j_r \neq n}} a_{j_1} \dots
a_{j_k} & a_n^{m_1} & \cdots & a_n^{m_{n-1}} \end{array} \right| \non \\
&=& \left| \begin{array}{cccc} a_2 \dots a_n & a_1^{m_1} & \cdots & a_1^{m_{n-1}} \\
\prod_{i \neq 2} a_i & a_2^{m_1} & \cdots & a_2^{m_{n-1}}  \\
\vdots & \vdots & \cdots & \vdots \\ a_1 \dots a_{n-1} & a_n^{m_1} & \cdots & a_n^{m_{n-1}} \end{array} \right| \non \\
&=& \sum_{j=1}^n (-1)^{j-1} \left| \begin{array}{ccc}  a_1^{m_1+1} & \cdots & a_1^{m_{n-1}+1} \\
\vdots & \cdots & \vdots \\  \big[\, a_j^{m_1+1}  & \cdots &  a_j^{m_{n-1}+1}\,\big]  \\ 
\vdots & \cdots & \vdots \\
a_n^{m-1+1} & \cdots & a_n^{m_{n-1}+1} \end{array} \right| = \left| \begin{array}{cccc}  1 & a_1^{m_1+1} & \cdots & a_1^{m_{n-1}+1} \\
\vdots & \vdots & \cdots & \vdots  \\ 1 & a_n^{m_1+1} & \cdots & a_n^{m_{n-1}+1} \end{array} \right|. \end{eqnarray*}
Next we continue by induction as before:
\begin{eqnarray*} \lefteqn{\left| \begin{array}{cccc} \sum_{\substack{1\leq j_1 < \dots < j_k \leq n \\ j_r \neq 1}}
a_{j_1} \dots a_{j_k} & a_1^{m_1} & \cdots & a_1^{m_{n-1}} \\ \sum_{\substack{1\leq j_1 < \dots < j_k \leq n \\ j_r
\neq 2}} a_{j_1} \dots a_{j_k} & a_2^{m_1} & \cdots & a_2^{m_{n-1}}  \\
\vdots & \vdots & \cdots & \vdots \\ \sum_{\substack{1\leq j_1 < \dots < j_k \leq n \\ j_r \neq n}} a_{j_1} \dots
a_{j_k} & a_n^{m_1} & \cdots & a_n^{m_{n-1}} \end{array} \right|} \\
&=& \sum_{j=1}^n (-1)^{j-1} \sum_{\substack{1\leq j_1 < \dots < j_k \leq n \\ j_r \neq j}} a_{j_1} \dots a_{j_k} \,
\left| \begin{array}{ccc}  a_1^{m_1} & \cdots & a_1^{m_{n-1}} \\
\vdots & \cdots & \vdots \\  \left[\, a_j^{m_1} \right. & \cdots & \left. a_j^{m_{n-1}}\, \right]  \\ \vdots & \cdots & \vdots \\
a_n^{m_1} & \cdots & a_n^{m_{n-1}} \end{array} \right| \\ 
&=& \sum_{j=1}^n (-1)^{j-1} \prod_{i\neq j} a_i^{m_1} \\ && \qquad \times \left|
\begin{array}{cccc} \sum_{\substack{1\leq j_1 < \dots < j_k \leq n \\ j_r \neq j}}
a_{j_1} \dots a_{j_k} & a_1^{m_2-m_1} & \cdots & a_1^{m_{n-1}-m_1} \\
\vdots & \cdots & \vdots \\  \bigg[\,\sum_{\substack{1\leq j_1 < \dots < j_k \leq n \\ j_r \neq j}}
a_{j_1} \dots a_{j_k} & a_j^{m_2-m_1} & \cdots &  a_j^{m_{n-1}-m_1}\bigg]  \\ 
\vdots & \vdots & \cdots & \vdots \\
\sum_{\substack{1\leq j_1 < \dots < j_k \leq n \\ j_r \neq j}} a_{j_1} \dots a_{j_k} 
& a_n^{m_2-m_1} & \cdots & a_n^{m_{n-1}-m_1} \end{array} \right| 
\end{eqnarray*}
Distinguishing the cases $j_r=i$ for some $r$ and $j_r\neq i$ for all $i$, this becomes
\begin{eqnarray*} && \sum_{j=1}^n (-1)^{j-1} \prod_{i\neq j} a_i^{m_1+1} \\ && \times \left|
\begin{array}{cccc} \sum_{\substack{1\leq j_1 < \dots < j_{k-1} \leq n \\ j_r \neq 1,j}}
a_{j_1} \dots a_{j_{k-1}}  & a_1^{m_2-m_1-1} & \cdots & a_1^{m_{n-1}-m_1-1} \\ \vdots & \cdots & \vdots \\
\left[\sum_{\substack{1\leq j_1 < \dots < j_{k-1} \leq n \\ j_r \neq j}}
a_{j_1} \dots a_{j_{k-1}} \right. &  a_j^{m_2-m_1-1} & \cdots & \left. a_j^{m_{n-1}-m_1-1}\right]  \\ \vdots & \vdots & \cdots & \vdots \\
\sum_{\substack{1\leq j_1 < \dots < j_{k-1} \leq n \\ j_r \neq j,n}} a_{j_1} \dots a_{j_{k-1}} &  a_n^{m_2-m_1-1} & \cdots & a_n^{m_{n-1}-m_1-1} \end{array} \right| \non \\ && 
+ \sum_{j=1}^n (-1)^{j-1} \prod_{i\neq j} a_i^{m_1} \\ && \times \left| \begin{array}{cccc}
\sum_{\substack{2\leq j_1 < \dots < j_k \leq n \\ j_r \neq j}} a_{j_1} \dots a_{j_k} & a_1^{m_2-m_1} & \cdots & a_1^{m_{n-1}-m_1} \\
\vdots & \cdots & \vdots \\  \bigg[\sum_{\substack{2\leq j_1 < \dots < j_k \leq n \\ j_r \neq j}} a_{j_1} \dots a_{j_k}  &
a_j^{m_2-m_1} & \cdots & a_j^{m_{n-1}-m_1}\bigg]  \\ \vdots & \vdots & \cdots & \vdots \\
\sum_{\substack{1\leq j_1 < \dots < j_k \leq n-1 \\ j_r \neq j}} a_{j_1} \dots a_{j_k} & a_n^{m_2-m_1} & \cdots & a_n^{m_{n-1}-m_1} \end{array} \right|. \end{eqnarray*} 
Both determinants in the last expression are of the same form as
the original, but of smaller size. By the induction hypothesis we therefore have that
\begin{eqnarray*} \lefteqn{\left| \begin{array}{cccc} 
\sum_{\substack{1\leq j_1 < \dots < j_k \leq n \\ j_r \neq 1}}
a_{j_1} \dots a_{j_k} & a_1^{m_1} & \cdots & a_1^{m_{n-1}} \\ 
\sum_{\substack{1\leq j_1 < \dots < j_k \leq n \\ j_r \neq 2}} 
a_{j_1} \dots a_{j_k} & a_2^{m_1} & \cdots & a_2^{m_{n-1}}  \\
\vdots & \vdots & \cdots & \vdots \\ \sum_{\substack{1\leq j_1 < \dots < j_k \leq n \\ j_r \neq n}} a_{j_1} \dots a_{j_k} & a_n^{m_1} & \cdots & a_n^{m_{n-1}} \end{array} \right|} \\
&=& \sum_{j=1}^n (-1)^{j-1} \prod_{i\neq j} a_i^{m_1+1} 
\sum_{\substack{m_2-m_1-1 \leq m'_1 < \dots < m'_{n-2}: \\
(\forall i) m'_i - m_{i+1} + m_1 = 0,-1 \\  \#\{i:\,m'_i = m_{i+1}-m_1\} = k-1}} \left| \begin{array}{cccc} 1 & a_1^{m'_1}  & \cdots & a_1^{m'_{n-2}} \\
\vdots & \vdots & \cdots & \vdots \\ \big[\,1  & a_j^{m'_1} & \cdots &  a_j^{m'_{n-2}} \big] \\ \vdots & \vdots & \cdots & \vdots \\
1 & a_n^{m'_1} & \cdots & a_n^{m'_{n-2}}  \end{array} \right| \non \\ && 
+ \sum_{j=1}^n (-1)^{j-1} \prod_{i\neq j}
a_i^{m_1} \sum_{\substack{m_2-m_1 \leq m'_1 < \dots < m'_{n-2}: \\
(\forall i) m'_i - m_{i+1} + m_1 = 0,1 \\  \#\{i:\,m'_i = m_{i+1}-m_1+1\} = k}} \left| \begin{array}{cccc} 1 & a_1^{m'_1}  & \cdots & a_1^{m'_{n-2}} \\
\vdots & \vdots & \cdots & \vdots \\ \big[\,1  & a_j^{m'_1} & \cdots & a_j^{m'_{n-2}} \big] \\ \vdots & \vdots & \cdots & \vdots \\
1 & a_n^{m'_1} & \cdots & a_n^{m'_{n-2}}  \end{array} \right| \\
&=& \sum_{j=1}^n (-1)^{j-1} \sum_{\substack{m_2 \leq m'_1 < \dots < m'_{n-2}: \\
(\forall i) m'_i - m_{i+1} = 0,1 \\  \#\{i:\,m'_i = m_{i+1}+1\} = k-1}} \left| \begin{array}{cccc} a_1^{m_1+1} & a_1^{m'_1}  & \cdots & a_1^{m'_{n-2}} \\
\vdots & \vdots & \cdots & \vdots \\ \big[\,a_1^{m_1+1} & a_j^{m'_1} & \cdots &  
a_j^{m'_{n-2}} \big] \\ \vdots & \vdots & \cdots & \vdots \\
a_n^{m_1+1} & a_n^{m'_1} & \cdots & a_n^{m'_{n-2}}  \end{array} \right| \non \\ 
&& + \sum_{j=1}^n (-1)^{j-1}
\sum_{\substack{m_2 \leq m'_1 < \dots < m'_{n-2}: \\
(\forall i) m'_i - m_{i+1} = 0,1 \\  \#\{i:\,m'_i = m_{i+1}+1\} = k}} 
\left| \begin{array}{cccc} a_1^{m_1} & a_1^{m'_1}  & \cdots & a_1^{m'_{n-2}} \\
\vdots & \vdots & \cdots & \vdots \\ \big[\,a_j^{m_1}  & a_j^{m'_1} & \cdots &  
a_j^{m'_{n-2}} \big] \\ \vdots & \vdots & \cdots & \vdots \\
a_n^{m_1} & a_n^{m'_1} & \cdots & a_n^{m'_{n-2}}  \end{array} \right|.
\end{eqnarray*}
Summing over $j$ we obtain
\begin{eqnarray*} \lefteqn{\left| \begin{array}{cccc} \sum_{\substack{1\leq j_1 < \dots < j_k \leq n \\ j_r \neq 1}}
a_{j_1} \dots a_{j_k} & a_1^{m_1} & \cdots & a_1^{m_{n-1}} \\ \sum_{\substack{1\leq j_1 < \dots < j_k \leq n \\ j_r
\neq 2}} a_{j_1} \dots a_{j_k} & a_2^{m_1} & \cdots & a_2^{m_{n-1}}  \\
\vdots & \vdots & \cdots & \vdots \\ \sum_{\substack{1\leq j_1 < \dots < j_k \leq n \\ j_r \neq n}} a_{j_1} \dots
a_{j_k} & a_n^{m_1} & \cdots & a_n^{m_{n-1}} \end{array} \right|} \\ &=&
\sum_{\substack{m_2 \leq m'_1 < \dots < m'_{n-2}: \\
(\forall i) m'_i - m_{i+1} = 0,1 \\  \#\{i:\,m'_i = m_{i+1}+1\} = k-1}} \left| \begin{array}{ccccc} 1 & a_1^{m_1+1} & a_1^{m'_1}  & \cdots & a_1^{m'_{n-2}} \\
\vdots & \vdots & \vdots & \cdots & \vdots \\ 1 & a_n^{m_1+1} & a_n^{m'_1} & \cdots & a_n^{m'_{n-2}}  \end{array}
\right| \\ && + \sum_{\substack{m_2 \leq m'_1 < \dots < m'_{n-2}: \\
(\forall i) m'_i - m_{i+1} = 0,1 \\  \#\{i:\,m'_i = m_{i+1}+1\} = k}} \left| \begin{array}{ccccc} 1 & a_1^{m_1} & a_1^{m'_1}  & \cdots & a_1^{m'_{n-2}} \\
\vdots & \vdots & \vdots & \cdots & \vdots \\ 1 & a_n^{m_1} & a_n^{m'_1} & \cdots & a_n^{m'_{n-2}}  \end{array} \right|
\\ &=& \sum_{\substack{m_1 \leq m'_1 < \dots < m'_{n-1}: \\
(\forall i) m'_i - m_{i} = 0,1 \\  \#\{i:\,m'_i = m_{i}+1\} = k}} \left| \begin{array}{cccc} 1 & a_1^{m'_1} & \cdots & a_1^{m'_{n-1}} \\
\vdots & \vdots & \cdots & \vdots \\ 1 & a_n^{m'_1} & \cdots & a_n^{m'_{n-1}}  \end{array} \right|.
\end{eqnarray*}
\qed

Iterating once more we get
\begin{cor} \label{L6cor} Let $\cal R$ be a commutative ring and $a_1,\dots,a_n \in {\cal R}$. Let $n \in \NN$ and $m_1,\dots,
m_{n-1} \in \NN$ such that $1 \leq m_1 < \dots < m_{n-1}$. Then, for any $k \in \NN$ with $1 \leq k \leq n$,
\begin{eqnarray} && \sum_{1 \leq j_1 < \dots < j_k \leq n} a_{j_1} \dots a_{j_k} \left| \begin{array}{cccc} 1 & a_1^{m_1} & \cdots & a_1^{m_{n-1}} \\
1 & a_2^{m_1} & \cdots & a_2^{m_{n-1}}  \\
\vdots & \vdots & \cdots & \vdots \\ 1 & a_n^{m_1} & \cdots & a_n^{m_{n-1}} \end{array} \right| \non \\
&& = \sum_{\substack{0 \leq m'_1 < \dots < m'_{n}: \, (\forall i) m'_i - m_{i-1} = 0,1 \\  \#\{i:\,m'_i = m_{i-1} + 1\} = k}}
\left| \begin{array}{ccc} a_1^{m'_1}  & \cdots & a_1^{m'_{n}} \\
\vdots & \cdots & \vdots \\ a_n^{m'_1} & \cdots & a_n^{m'_{n}}  \end{array} \right|,
\end{eqnarray}
where $m_0=0$.
\end{cor}


\end{document}